\documentclass[11pt]{amsart}

\usepackage[square,compress,comma, numbers,sort]{natbib}
\usepackage[colorlinks=true, citecolor=red, linkcolor=blue]{hyperref}
\usepackage{amsfonts,mathtools,mathabx}
\usepackage[shortlabels]{enumitem}
\usepackage{cleveref}
\allowdisplaybreaks[4]

\def\cadlag{c\`adl\`ag }

\usepackage{amssymb}
\usepackage{color}
\newcommand{\abs}[1]{\left\lvert #1 \right\rvert}

\usepackage{amsfonts,mathtools}

\DeclarePairedDelimiterXPP\pk[1]{\mathbb{P}}\{ \}{}{ #1}
\DeclarePairedDelimiterXPP\E[1]{\mathbb{E}}\{ \}{}{	#1}

\def\FRE{\mbox{Fr\'{e}chet }}

 \usepackage{xparse}% http://ctan.org/pkg/xparse

\NewDocumentCommand{\ceil}{s O{} m}{%
  \IfBooleanTF{#1} % starred
    {\left\lceil#3\right\rceil} % \ceil*[..]{..}
    {#2\lceil#3#2\rceil} % \ceil[..]{..}
}
\NewDocumentCommand{\floor}{s O{} m}{%
  \IfBooleanTF{#1} % starred
    {\left\lfloor#3\right\rfloor}
    {#2\lfloor#3#2\rfloor}
}

\newcommand{\norm}[1]{\lVert #1 \rVert}

\definecolor{c20}{rgb}{0.,0.7,0.}
\definecolor{c30}{rgb}{0.,0.,1.}
\definecolor{c40}{rgb}{1,0.1,0.7}
\definecolor{c50}{rgb}{1,0,0}
\definecolor{c60}{rgb}{1,0.9,0.1}
\definecolor{c70}{rgb}{0.50,1.00,0.00}
\def\cEH#1{{\textcolor{c30}{#1}}}
\def\cEH#1{#1}

\def\bE#1{{\textcolor{c30}{#1}}}

\def\cE#1{{\textcolor{c30}{#1}}}
\def\cE#1{#1}

\def\zm#1{#1}

\def\kdd#1{{\textcolor{red}{#1}}}
\definecolor{orange}{rgb}{1,0.5,0}
\def\k2#1{{\textcolor{orange}{#1}}}
\numberwithin{equation}{section}
\newtheorem{theo}{Theorem}[section]
\newtheorem{sat}[theo]{Proposition}
\newtheorem{de}[theo]{Definition}
\newtheorem{lem}[theo]{Lemma}

\newtheorem{example}[theo]{Example}
\newtheorem{korr}[theo]{Corollary}
\newtheorem{remark}[theo]{Remark}

\numberwithin{equation}{section}

\newcommand{\prooftheo}[1]{ \textsc{Proof of Theorem} \ref{#1} }
\newcommand{\proofprop}[1]{\textsc{Proof of Proposition} \ref{#1}}
\newcommand{\prooflem}[1]{\textsc{Proof of Lemma} \ref{#1}}
\newcommand{\proofkorr}[1]{\textsc{Proof of Corollary} \ref{#1}}

\newcommand{\QED}{\hfill $\Box$}

\newcommand{\COM}[1]{}

\def\IF{\infty}

\newcommand{\R}{\mathbb{R}}
\newcommand{\inr}{\in \R}

\newcommand{\BQN}{\begin{eqnarray}}
\newcommand{\EQN}{\end{eqnarray}}
\newcommand{\BQNY}{\begin{eqnarray*}}
\newcommand{\EQNY}{\end{eqnarray*}}

\newcommand{\limit}[1]{\lim_{#1 \to   \infty}}

\def\bqny#1{\begin{eqnarray*} #1 \end{eqnarray*}}
\def\bqn#1{\begin{eqnarray} #1 \end{eqnarray}}

\newcommand{\BS}{\begin{sat}}
\newcommand{\ES}{\end{sat}}
\newcommand{\BT}{\begin{theo}}
\newcommand{\ET}{\end{theo}}
\newcommand{\BK}{\begin{korr}}
\newcommand{\EK}{\end{korr}}

\newcommand{\BEX}{\begin{example}}
\newcommand{\EEX}{\end{example}}

\newcommand{\BD}{\begin{de}}
\newcommand{\ED}{\end{de}}

\newcommand{\BIT}{\begin{itemize}}
\newcommand{\EIT}{\end{itemize}}
\newcommand{\BDI}{\begin{description}}
\newcommand{\EDI}{\end{description}}

\newcommand{\BRM}{\begin{remark}}
\newcommand{\ERM}{\end{remark}}

\newcommand{\BEL}{\begin{lem}}
\newcommand{\EEL}{\end{lem}}

\newcommand{\nelem}[1]{\Cref{#1}}

\newcommand{\netheo}[1]{\Cref{#1}}

\newcommand{\nekorr}[1]{\Cref{#1}}

\def\TT{\mathcal{T}}
\def\TT{\R }

\def\Z{\mathbb{Z}}
\def\inn{\in \mathbb{N}}

\newcommand\ind[1]{\mathbb{I}{\left\{#1\right\}}}

\def\TT{\R}
\def\intT{\int_{\TT}}

\begin{document}

\title[Berman functions]
{On Berman functions}

\author{Krzysztof D\c{e}bicki}
\address{Krzysztof D\c{e}bicki, Mathematical Institute, University of Wroc\l aw, pl. Grunwaldzki 2/4, 50-384 Wroc\l aw, Poland}
\email{Krzysztof.Debicki@math.uni.wroc.pl}

\author{Enkelejd  Hashorva}
\address{Enkelejd Hashorva, Department of Actuarial Science, University of Lausanne,\\
Chamberonne 1015 Lausanne, Switzerland}
\email{Enkelejd.Hashorva@unil.ch}

\author{Zbigniew Michna}
\address{Zbigniew Michna, Department of Operations Research and Business Intelligence, Wrocław University of Science and Technology, Poland}
\email{Zbigniew.Michna@pwr.edu.pl}

\bigskip

\date{\today}
 \maketitle

\begin{quote}

{\bf Abstract:}
Let  $Z(t)= \exp\left( \sqrt{ 2} B_H(t)- \abs{t}^{2H}\right), t\inr$ with $B_H(t),t\inr$  a
standard fractional Brownian motion (fBm) with Hurst parameter  $H \in (0,1]$ and define for $x$ non-negative the Berman function
\bqny{
	\mathcal{B}_{Z}(x)= %= \int_x^\infty\frac{1}{y}dG(y)=
	\E*{ \frac{ \ind{ \epsilon_0(RZ) > x}}{ \epsilon_0(RZ)}}\in (0,\IF),
}
where the random variable   $R$ independent of $Z$ has  survival function $1/x,x\ge 1$
 and
$$
\epsilon_0(RZ) = \int_{\R} \ind{ RZ(t)> 1}\zm{dt} .
$$
 In this paper we consider a general random field (rf) $Z$ that is a spectral rf of some stationary max-stable rf $X$ and derive the properties of the corresponding
Berman functions.  In particular, we show that Berman functions can be approximated by the corresponding discrete ones and derive interesting representations of those functions which are of interest for Monte Carlo simulations, which are presented in this article.
\end{quote}

{\bf Key Words:} Berman functions; Pickands constants; max-stable random fields;  simulations

{\bf AMS Classification:} Primary 60G15; secondary 60G70\\

\section{Introduction}
In the study of sojourns of rf's in a series of papers by S.\ Berman, see e.g., \cite{Berman82,Berman92} a key random variable (rv) and a  related constant appear.
Specifically, let
$Z(t)= \exp( \sqrt{ 2} B_{\cE{H}}(t)- \abs{t}^{2 H}), t\inr,$ with $B_{H}$ a
fractional Brownian motion (fBm) with Hurst parameter $H\in(0,1]$,
that is a centered Gaussian process \cE{with stationary increments,  $Var(B_H(t))=\abs{t}^{2H},t\inr$ and continuous sample paths.}
In view of \cite[Thm 3.3.1, Eq.\ (3.3.6)]{Berman92} the following rv (hereafter 
 $\ind{\cdot}$ is the indicator function)
$$
\epsilon_0(RZ) = \int_{\R} \ind{ RZ(t)> 1}      dt
$$
 plays a crucial role  in the analysis of extremes of Gaussian processes. Throughout this paper $R$ is a 1-Pareto rv ($\ln R$ is unit exponential) independent of any other random element.\\
The distribution function of $\epsilon_0(RZ)$ is known only for  {$H\in \{1/2,1\}$}.  For  {$H=1$} as shown in \cite[Eq.\ (3.3.23)]{Berman92} $\epsilon_0(RZ)$
has probability density function (pdf) $ x^2 e^{- x^2/2}/(2 \sqrt{ \pi}),x>0$, whereas for {$H=1/2$} its  pdf is calculated in \cite[Eq.\ (5.6.9)]{Berman92}.\\
The so-called Berman function   defined for all $x\ge 0$  (see \cite[Eq.\ (3.0.2)]{Berman92})  given by
\bqn{\label{bermanrep}
	\mathcal{B}_{Z}(x) %= \int_x^\infty\frac{1}{y}dG(y)
	&=&
	\E*{ \frac{ \ind{ \epsilon_0(RZ) > x}}{ \epsilon_0(RZ)}}\in (0,\IF)
}
%where
%$$
%G(y) =\pk*{ \epsilon_0(RZ) \leq y}
%$$
is defined in \cite[Thm 3.3.1, Eq.\ (3.3.6)]{Berman92}.

An important property of the Berman function is that for $x=0$ it equals the {\it Pickands} constant, see \cite[Thm 10.5.1]{Berman92} i.e.,
$\mathcal{B}_{Z }(0) = \mathcal{H}_{Z }, $ where $\mathcal{H}_{Z } $ is the so called {\it generalised Pickands} constant
\bqny{
\mathcal{H}_{Z } =
\limit{T} \frac{1}{T} \E*{  \sup_{ t\in [0,T]} Z (t)}.
}
 This fact is crucial since $\mathcal{B}_{Z }(0)$  is the first known expression of $\mathcal{H}_{Z } $
 in terms of an expectation, which is of particular usefulness for simulation purposes, see \cite{Faletal2010,DiekerY,MR3679987} for details on classical Pickands constants.\\
  Besides,  Berman's representation of Pickands constant yields tight  lower bounds for $\mathcal{H}_{Z } $, see 
 \cite[Thm 1.1]{KSojourn}. As shown in  \cite{KSojourn} for all $x\ge 0$
$$
\mathcal{B}_{Z }(x)
= \limit{T} \frac{1}{T} \mathcal{B} ([0,T], x) , \quad
\mathcal{B}_{Z } ([0,T], x ) \coloneqq
\int_0^\IF  \pk*{  \int_{0}^T \mathbb{I}(Z (t) >   s) dt > x }ds.
$$ 
Motivated by the above definition, in this contribution we shall introduce the Berman functions for given $\delta\ge 0$
  with respect to some non-negative rf
  $Z(t),t\in \R^d,d\ge 1$   with \cadlag  sample paths
(see e.g., \cite{Svante,MartinE}
  for the definition and properties of
generalised \cadlag functions)
such that 
\bqn{\E{ Z (t) }=1,  \quad t\in \R^d.
	\label{MS}
}
Specifically, for given  non-negative $\delta, x$ define
$$
\mathcal{B}_{Z  }^\delta(x)\coloneqq\lim_{T\rightarrow \infty}\frac{1}{T^{\zm d}}
\mathcal{B}_{Z  }^\delta([0,T]^d \cap \delta \Z^d, x),
$$
where
$$
\mathcal{B}_{Z  }^\delta([0,T]^d \cap \delta \Z^d, x)\coloneqq
 \int_0^\IF \pk*{  \int_{[0,T]^d \cap \delta \Z^d} \mathbb{I}(Z (t) >   s) \lambda_\delta(dt) > x }ds.
$$
Here 
$\lambda_0(dt)=\lambda(dt)$ is the Lebesgue measure on $\R^d$, $0 \Z^d=\R^d$  and  $\lambda_\delta(dt)/\delta^d$
is the counting measure on $\delta \Z^d$ if $\delta >0$. Hence $\mathcal{B}_{Z  }^\delta(x),\delta>0$ is the discrete counterpart of
$\mathcal{B}_{Z  }(x)$ and $\mathcal{B}_{Z  }^0(x)=\mathcal{B}_{Z  }(x)$.

In general, in order to be well-defined for the function  
$\mathcal{B}_{Z}^\delta(x), x\ge 0$ some further restriction on the rf $Z$ are needed.
A very tractable case for which we can  utilise results from the theory of max-stable stationary rf's is when  $Z$ is the spectral rf of a stationary max-stable rf $X(t), t\in \R^d$, see \eqref{eq1} below.\\
An interesting special case is when $\ln Z(t)$ is a Gaussian rf with trend equal the half of its variance function having further stationary increments. We shall show in  \Cref{l.sojourn} that for such $Z$
the corresponding Berman function $\mathcal{B}_{Z  }(x)$ appears in the tail asymptotic of the sojourn of a related Gaussian rf.

Organisation of the rest of the paper.
In Section 2  we first present in \Cref{erdha} a formula for Berman functions
and then in \Cref{K1} and \Cref{korrDe} we show some continuity properties of those functions.
 In \Cref{repBE} and \Cref{xhxh} we present two representations for Berman functions and discuss conditions for their positivity. Section 3 is dedicated to the approximation of Berman functions focusing on the Gaussian case. All the proofs are postponed to Section 4.

\def\HDM{ \mathcal{H}_{X,r}^\delta}

\def\HDO{ \mathcal{H}_{X,r}}

\def\TT{\mathcal{T}}
\def\TT{\R^d}
\def\AA{\mathcal{D}}

\section{Main Results}
Let the rf $Z(t), t\in \TT$ be as above  defined in the
 non-atomic complete probability space  $(\Omega, \mathcal{F}, \mathbb{P})$.
 Let further $X(t),t\in \TT $ %with $\TT= \R^d$
be a max-stable stationary rf, which has spectral rf  $Z$ in its de Haan representation  (see e.g., \cite{deHaan,booksoulier})
\bqn{\label{eq1}
	X(t)=  \max_{i\ge 1} \Gamma_i^{{-1}}  Z^{(i)}(t), \quad t\in \TT.
}
Here $\Gamma_i= \sum_{k=1}^i \mathcal{V}_k$ with $\mathcal{V}_k, k\ge 1$ mutually independent
unit exponential rv's being independent of $\{Z^{(i)}\}_{i=1}^\infty$ which are independent copies of $Z$. For simplicity we shall assume that the marginal distributions of the rf $X$ are unit \FRE (equal to $e^{-1/x },x>0$) which in turn implies $\E{Z (t)} =1$ for all $t\in \TT$.

Suppose further that for all $T>0$
\bqn{\label{norm}
\E*{ \sup_{t\in [0,T]^d}  Z (t)}< \IF
}
and $Z$ has almost surely sample paths on the space $D$ of non-negative
\cadlag functions $f: \TT \mapsto [0, \IF)$ equipped with Skorohod's $J_1$-topology. We  shall denote by
$\mathcal{D}=\sigma(\pi_t, t\in T_0)$ the $\sigma$-field generated by the projection maps $\pi_t: \pi_t f= f(t), f\in D$ with $T_0$ a countable dense subset of $\R^d$.
In view of  \cite[Thm 6.9]{Htilt} with $\alpha=1, L=B^{-1}$, see also \cite[Eq.\ (5.2)]{Hrovje} the stationarity of  $X$ is equivalent with
\bqn{
	\label{rinashero}
	 \E{ Z (h) F(Z)} = \E{ Z (0) F(B^h Z)}, \quad \forall h\in \R^d
}
valid for every  measurable functional $F: D \to [0,\IF]$ such that $F(cf)= F(f)$ for all $f\in D, c>0$. Here  we use the standard notation $B^h Z(\cdot )= Z(\cdot-h), h\inr^d$.

We shall suppose next without loss of generality (see \cite[Lem 7.1]{HBernulli}) that
\bqn{ \pk*{\sup_{t\in \R^d}  Z(t) >0}=1.
	\label{sid}
}
 Under the assumption that $X$ is stationary $\mathcal{B}_{Z }^\delta(x)$  is well-defined for all $\delta,x$ non-negative as we shall show below.  We note first that, see e.g., \cite{KSojourn,ZKE}
$$
 \limit{T} {\frac{1}{T^{d}}}\mathcal{B}_{Z}^\delta([0,T]^{d}\cap \delta \Z^{\zm{\delta}} ,0)=\mathcal{B}_{Z}^\delta(0) =\mathcal{H}_{Z}^\delta \in (0, \IF),
$$
 where $\mathcal{H}_{Z}^\delta$ is the discrete counterpart of the classical Pickands constant  $\mathcal{H}_{Z}=\mathcal{H}_{Z}^0$. Hence  for any $x>0$ we have
 $$\mathcal{B}_{Z}^\delta(x) \le \mathcal{H}_{ Z }^\delta< \IF. $$
Set below for $\delta>0$
$$ S_\delta=S_\delta(Z)=\int_{ \delta \Z^d}  Z (t) \lambda_\delta(dt)=
\delta^d \sum_{ t\in \delta \Z^d}  Z (t) 
$$
and let $S_0=S_0(Z)=\int_{ \R^d}  Z (t) \lambda(dt)$. 
 In view of \eqref{sid} we have that $S_0>0$ almost surely.
Since we do not consider the case $\delta>0$ and $\delta=0$ simultaneously, we can assume that $S_\delta>0$ almost surely (we can construct a spectral rf $Z$ for $X$ that guarantees this, see \cite[Lem 7.3]{HBernulli}).

In view of \cite[Cor 2.1]{ZKE} if  $\pk{ S_0=\IF} =1$, then  $\mathcal{H}_{ Z }=0$ implying 
 $$\mathcal{B}_{ Z }^\delta(x){=}\mathcal{H}_{ Z  }=0, \quad \forall \delta, x\ge 0.$$
 The next result states the existence and the positivity of Berman functions presenting 
 further  a tractable formula that is useful for simulations of those functions.
\BT    If $\pk{ S_0=\IF} <1$, then for any $\delta,   x$ non-negative constants we have
\bqn{ \label{bli}
	\mathcal{B}_{ Z }^\delta(x) &=& \int_0^\IF \E*{ \frac{ Z (0)}{ S_{\cE{\delta}}} \mathbb{I} \Bigl(  \int_{\delta \Z^d}  \mathbb{I}(  Z   (t) >  s)  \lambda_\delta(dt) > x \Bigr) }\lambda(ds) < \IF.
}
 Moreover,  \eqref{bli} holds substituting $S_\delta$ by $S_\eta$,
  where $\eta>0$ if $\delta=0$ and $\eta=k \delta, k\inn$ if $\delta>0$, provided that
\bqn{ \label{laps}
	\{ {S_0(Z)} < \IF \} \subset   \{ S_\eta(B^r Z) \in (0,\IF)   \}, \quad  \forall r\in \delta \Z^d
}
almost surely.
\label{erdha}
\ET

\BRM
\begin{enumerate}[(i)]
	\item If $x=0$, then we retrieve the results of \cite[Prop 2.1]{ZKE}.
	\item As shown in \cite{ZKE} condition \eqref{laps} holds in the particular case that $Z(t)>0,t\inr^d$ almost surely.
	\item  \label{MMS} One example for $Z$, see for instance \cite{ZKE} is taking
	$$Z(t)= \exp( V(t)- \sigma^2_V(t)/2)), \quad t\inr^d,$$
	 where $V(t),t\inr^d$ is a centered Gaussian rf with almost surely continuous trajectories and stationary increments, $\sigma^2_V(t)=Var(V(t))$ and $\sigma_V(0)=0$. For this case $Z(t)>0,t\inr^d$ almost surely,  condition  \eqref{laps} is satisfied and \eqref{bli} reads
	\bqn{ \label{bli2}
		\mathcal{B}_{ Z }^\delta(x) &=& \int_0^\IF \E*{ \frac{1}{ S_{\cE{\delta}}} \mathbb{I} \Bigl(  \int_{\delta \Z^d}  \mathbb{I}(  V(t) - \sigma^2_V(t)/2>  \ln s)  \lambda_\delta(dt) > x \Bigr) }\lambda(ds) < \IF.
	}		
\end{enumerate}
\label{cern}
\ERM

\BK If $Z$ has almost surely continuous trajectories, then for  all  $x_0\geq 0$
\bqn{
	\lim_{ x\rightarrow x_0} 	\mathcal{ B}_{ Z }^ 0(x)= 	\mathcal{ B}_{ Z }^ 0(x_0).
}
\label{K1}
\EK

Define next a probability measure $\mu$ on $\AA$ by
\bqn{ \label{thet}
	\mu(A)= \E{  Z (0)\ind{ Z/Z(0) \in A}}, \quad A \in \AA.
}	
Let $\Theta$ be a rf with law $\mu$. By the definition, $\Theta$ has also \cadlag\ sample paths
and since $D$ is  Polish,  in view of \cite[Lem p.\ 1276]{MR100291} we can assume that $\Theta$ is defined in the same probability space as $Z$. Recall that $\lambda_\delta(dt)/\delta^d$
is the counting measure on $\delta \Z^d$ if $\delta >0$ and $\lambda_0$ is the Lebesgue measure on $\R^d$.
 Since we can rewrite \eqref{bli} as
\bqn{\label{eqTE}
	\mathcal{ B}_{ Z }^ \delta(x) &=&
	 \int_0^\IF\E*{ \frac{1}{ S_{\cE{\delta}}(\Theta) }  \mathbb{I} \Bigl(  \int_{\delta \Z^d}  \mathbb{I}( \zm{\Theta}  (t) >  s)  \lambda_\delta(dt) > x \Bigr)  }\lambda(ds) < \IF,
}
where
$$S_{{\delta}}(\Theta)= \int_{\cE{\delta  \Z^d}} \Theta(t) \lambda_\delta(dt)$$
 and the law of $\Theta$ is uniquely determined by the law of the max-stable stationary rf $X$ and does not depend on the particular choice of $Z$, see  \cite[Lem A.1]{Htilt},
\cE{hence if $Z_*$} is another spectral rf for $X$, then
\bqn{\label{eqber}
	\mathcal{ B}_{ Z }^ \delta(x)= \mathcal{ B}_\zm{{Z_*}}^ \delta(x)
}	
for all  $\delta\geq 0$. Assume next that $\pk{S_0< \IF}=1$  and  let
\bqn{ Z_*(t)= (p(T))^{-1}B^T  {\bar Q_\delta(t)}, \quad t\in \delta \Z^d
\label{Q}
}
be  a spectral rf of  the max-stable rf $X_\delta(t)= X(t), t\in \delta \Z^d$,  where \cE{$\bar Q_\delta$ is independent of a rv $T$}, which has pdf $p(s)>0, \,s\in \cE{\delta \Z^d}$. We choose $p$ to be continuous when $\delta=0$. In view of \cite[Thm \kdd{2.3}]{debicki2017approximation} one possible construction
  is
  $$\bar Q_\delta(t)=   c\frac{\Theta(t)}{  S_\delta(\Theta)} , \quad t\in \delta \Z^d,$$
   \cE{with $c=1$ if $\delta=0$ and $c=\delta^d$ otherwise}.
\cE{Set below $Q_\delta= \bar Q_\delta/c$.}
\BEL \begin{enumerate}[(i)]
	\item \label{en:A} If $\pk{S_0< \IF}=1$, \cE{then for  $  Q_\delta$ as above and all}  $\delta,x$ non-negative we have
\bqn{
	\mathcal{ B}_{ Z }^ \delta(x) &=&  \int_0^\IF \pk*{  \int_{\delta \Z^d}  \mathbb{I}( Q_\delta  (t) >  s)  \lambda_\delta(dt) > x }\lambda(ds)  < \IF.
}
\item \label{en:B} If $\pk{ S_0< \IF}>0$, then with  $V(t)= Z(t) \lvert  S_0 < \IF $
for all $\delta,x$ non-negative we have
\bqn{
\label{paP}
	\mathcal{ B}_{ Z }^ \delta(x) &=& \pk{S_0(\Theta)< \IF  } \mathcal{ B}_{V}^ \delta(x) < \IF.
}
\end{enumerate}
\label{xhxh}
\EEL
Let in the following $Y(t)= R \Theta(t)$   with $R$   a $1$-Pareto rv with survival function $1/x, x \ge 1$ independent of $\Theta$ and set hereafter 
$$
\epsilon_\delta(Y)=  \int_{\delta \Z^d} \ind{ Y(t)> 1} \lambda_\delta(dt).
$$
Recall that when $\delta=0$ we interpret $\delta \Z^d$ as $\R^d$. We establish below the Berman representation \eqref{bermanrep} for the general setup of this paper.
\BT\label{phl}
 If $\pk{ S_0=\IF} <1$, then  for all $ \delta,  x$ non-negative
\bqn{ 	\mathcal{ B}_{ Z }^\delta (x) =   \E*{\frac{ \ind{ \epsilon_\delta(Y )> x }}{\epsilon_\delta (Y )}      }  < \IF.
	\label{repBE}
}
\ET

\BK \label{korrIm} Under the conditions of \netheo{phl} we have that $\epsilon_0 (Y )$ \cE{has a continuous distribution} if $Z$ \cE{has almost surely continuous trajectories}. Moreover,  $\mathcal{ B}_{ Z }^\delta (x)    >0$  for all $x\ge 0$ such that $\pk{ \epsilon_\delta (Y ) >x} >0$.
\EK

\begin{sat}\label{ineq1} For all  $\delta\ge 0$ and $x>0$  we have
	\bqn{\label{eqXhK}
 		\frac{\pk{\epsilon_\delta(Y)>x}^2}{\E{\epsilon_\delta(Y)}} \le \mathcal{ B}_{ Z }^\delta (x)  \le x^{-1}\pk{\epsilon_\delta(Y)\geq x}   .
	}
%	and
%	\bqn{\label{upp}
%		\mathcal{ B}_{ Z }^\delta (x) \leq \frac{}{x}.
%	}	
\end{sat}
\BRM\label{r.42}
\begin{enumerate}[(i)]
	\item
 If $x=0$, the lower bound in \eqref{eqXhK} holds with 1 in the numerator, see \cite{KSojourn,HBernulli}.
\item If $\E{\epsilon^p_\delta(Y)}$ is finite for some $p>0$, then combination of the upper bound in \eqref{eqXhK} with the Markov inequality
gives the following upper bound 
\bqn{\label{upBMar}
	\mathcal{ B}_{ Z }^\delta (x) \leq x^{-p-1}\E{\epsilon^p_\delta(Y)}, \quad x>0.
}
\item If $\E{\epsilon_\delta(Y)}<\infty$ and $\int_0^\infty e^{sx}(\mathcal{ B}_{ Z }^\delta (x))^{1/2}dx<\infty$, then it follows that for all $s>0$
$$
\E{e^{s\epsilon_\delta(Y)}}\leq 1+s(\E{\epsilon_\delta(Y)})^{1/2}\int_0^\infty e^{sx}(\mathcal{ B}_{ Z }^\delta (x))^{1/2}dx.
$$
\item \label{ZV} Since $Y= R\Theta$, we can calculate in case of known $\Theta$ the expectation of $\epsilon_\delta(Y)$ as follows
\begin{eqnarray*}
	\E{\epsilon_\delta(Y)}%&=&\int_{\delta \Z^d} \E{\ind{ Y(t)> 1}} \lambda_\delta(dt)\\
	%&=&\int_{\delta \Z^d} \pk{ Y(t)> 1} \lambda_\delta(dt)\\
	&=&\int_{\delta \Z^d} \pk{ R\Theta(t)> 1} \lambda_\delta(dt)
	%&=&\int_{\delta \Z^d} \int_1^\infty\pk{ \Theta(t)> 1/r}r^{-2}dr \lambda_\delta(dt)\\
	=\int_1^\infty \int_{\delta \Z^d}\pk{ \Theta(t)> 1/r}\lambda_\delta(dt)r^{-2}dr.
\end{eqnarray*}
If $Z(t)= \exp( V(t)- \sigma^2_V(t)/2)), t\inr^d$ is as in \Cref{cern}, \Cref{MMS},
% where $V(t),t\inr^d$ is a centered Gaussian rf with almost surely
%continuous trajectories and stationary increments, $\sigma^2_V(t)=Var(V(t))$ and $\sigma_V(0)=0$.
then in view of \cite[Lem 5.4]{KSojourn}, \cite[Eq.\ (5.3)]{HBernulli} we have
\bqn{  \E{\epsilon_\delta(Y)}=  \int_{\delta \Z^d}\int_1^\infty {\Psi}\left(\frac{\sigma_V(t)}{2}-\frac{\ln  r}{\sigma(t)}\right)r^{-2}\lambda_\delta(dt)dr = 2  \int_{t\in \delta \Z^d} {\Psi}(\sigma_V(t)
	/2)\lambda_\delta(dt),
	\label{expepy}
}
where $\Psi$ is the survival function of an $N(0,1)$ rv.
\end{enumerate}
\ERM

%Consider next the log-Gaussian case as in \Cref{s.rate}, i.e.,

%In particular,
%for $p=1$

\section{Approximation of $\mathcal{B}_{Z }^\delta(x)$ and its behaviour for large $x$}
\label{s.rate}
%In this section we discuss the approximation of $\mathcal{ B}_{ Z }^\delta$.
We show first that $\mathcal{ B}_{ Z }=\mathcal{ B}_{ Z }^{0}$ can be approximated by considering $\mathcal{ B}_{ Z }^{ \delta }(x)$ and letting $\delta \downarrow 0$.
\BS For all $x\ge 0$ we have that
$$ \lim_{\delta \downarrow 0}  \mathcal{ B}_{ Z }^{ \delta } (x) =  \mathcal{ B}_{ Z }^{0}  (x) .
$$
\label{korrDe}
\ES
We note in passing that for $x=0$ we retrieve the approximation for Pickands constants derived in \cite{ZKE}.
 An approximation of $\mathcal{ B}_{ Z }^{ \delta }(x)$ can be obtained by letting $T\to \IF$
and calculating the limit of
\[\frac{\mathcal{B}^\delta_{Z }([0,T]^d\cap \delta \Z^d,x)}{ T^d}.
\]
For such an approximation we shall discuss the rate of convergence to $\mathcal{ B}_{ Z }^{ \delta }(x)$ assuming further
that
$$
Z(t)=\exp\left(  V(t)- \frac{\sigma^2_V(t)}{2}\right), \quad  t\in\TT
$$
is as in \Cref{cern}, \Cref{MMS}.
% Hereafter  $V$ is an almost surely continuous
%centered Gaussian field with stationary increments, $V(0)=0$ a.s.\ and variance function $\sigma^2_V(t)\coloneqq Var(V(t))$
%that satisfies the following conditions
%\\

{\bf A1} $\sigma^2_V(t)$ is a continuous and strictly increasing function, and
there exists $\alpha_0\in(0,2]$ and $A_0\in (0,\infty)$ such that
\[
\limsup_{\lVert t\rVert\to0}\frac{\sigma^2_V(t)}{\lVert t\rVert^{\alpha_0}}\le A_0,
\]
where $\lVert \cdot\rVert$ is the Euclidean norm.
\\
{\bf A2} There exists $\alpha_\infty\in (0,2]$ such that
\[
\liminf_{\lVert t\rVert\to\infty} \frac{\sigma^2_V(t)}{\lVert t\rVert^{\alpha_\infty}}>0.
\]
The following theorem constitutes the main finding of this section.
\BT\label{Th.rate}
Under {\bf A1-A2} we have  for all $\delta, x$ non-negative and $\lambda \in (0,1)$
\begin{eqnarray}
\lim_{T\to\infty}\left|\mathcal{ B}_{ Z }^\delta(x)-\frac{\mathcal{ B}_{ Z }^\delta([0,T]^d\cap \delta \Z^d,x)}{T^d}\right|
T^\lambda=0.
\label{eq:rate}
\end{eqnarray}
\ET

\BRM
\begin{enumerate}[(i)]
	\item For $x=0$ the rate of convergence in \eqref{eq:rate} agrees with the findings in \cite{MR2222683}.
\item The range of the parameter $\lambda\in(0,1)$ in Theorem \ref{Th.rate} cannot be extended to $\lambda\ge1$.
Indeed, following \cite{Ch19}, for $V(t)=\sqrt{2}B_{1}(t)$, $\delta=0$, $T>x$ and $d=1$ we have
\[
\mathcal{ B}_{ Z }([0,T],x)=
2\Psi(x/\sqrt{2})+\sqrt{2}(T-x)\varphi(x/\sqrt{2})
\]
implying
\BQN\label{leads}
\mathcal{ B}_{ Z }(x)=\sqrt{2}\varphi(x/\sqrt{2}),
\EQN
where $\varphi(\cdot)$ is the pdf of an $N(0,1)$ rv. Consequently, we have
\begin{eqnarray*}
\lim_{T\to\infty}
\left|
\mathcal{ B}_{ Z }(x)-\frac{\mathcal{ B}_{ Z }([0,T],x)}{T}
\right|T
={ 2\Psi(x/\sqrt{2}) } >0. %-\sqrt{2}x \varphi(x/\sqrt{2})
\end{eqnarray*}
\end{enumerate}
\ERM

%and $f:\R^+\to\R^+$ be such that
%$f(x)<x$ for each $x\ge 1$.
%If ${\Lambda}>0$ is such that
%$
%\sup_{t> 0}\frac{\sigma_{V}^2(t)}{\sigma_{V}^2({\Lambda}t)}\le \frac{1}{2},
%$
%then for each $T\in \mathbb{N}$ and each $x\ge0$
%\begin{eqnarray}
%\left|\mathcal{ B}_{ Z }(x)-\frac{\mathcal{ B}_{ Z }([0,T],x)}{T}\right|&\le&
%\frac{\mathcal{H}_{Z}( f\tilde{}(T))}{T}
%+
%\frac{\mathcal{H}_{Z}^2(\Lambda^5T)}{T}
%\exp\left(-\frac{\sigma^2_{V}( f(T))}{8}\right)
%+\nonumber\\
%& &+
%\frac{\mathcal{H}^2_{Z}(\Lambda^5T)}{T}\sum_{i=1}^\infty
%\exp\left(-\frac{\sigma^2_{V}(iT)}{8}\right).
%\end{eqnarray}

%\section{Estimates of $\mathcal{ B}_{ Z }^\delta (x)$}

%In this section we investigate some properties of the Berman function
%$\mathcal{ B}_{ Z }^\delta (x)$.
%Following Theorem \ref{phl} we begin with an observation that
%$\mathcal{ B}_{ Z }^\delta (x)$ is intimately related with
%$\pk{\epsilon_\delta(Y)>x}$.

In the rest of this section we focus on $d=1$ log-Gaussian case.
%	Recall that
%for $V(t)=\sqrt{2}B_{1/2}(t)$ we have
%\begin{eqnarray}
%\mathcal{ B}_{ Z } (x)=(2+x)\Psi\left(\sqrt{\frac{x}{2}}\right)-\sqrt{\frac{x}{\pi}}\exp\left(-\frac{x}{4}\right)\label{b.bm}
%\end{eqnarray}
%and for
%$V(t)=\sqrt{2}B_{1}(t)$
%Berman function reads
%\begin{eqnarray}
%\mathcal{ B}_{ Z } (x)=\sqrt{2}\varphi(x/\sqrt{2}). \label{b.dbm}
%\end{eqnarray}
In view of \eqref{leads} for some finite positive constant $C$
\[
\ln  (\mathcal{ B}_{ Z }^\delta (x)) \sim -{ C} \sigma^2_V(x), \quad x\to \IF.
\]
The next result gives logarithmic bounds for $\mathcal{ B}_{ Z }^\delta (x)$
as $x\to\infty$ that supports this hypothesis.
\begin{sat}\label{prop.bound}
Suppose that $d=1$ and $V$ satisfies {\bf A1-A2}.  Then
\[
\liminf_{x\to\infty}\frac{\ln (\mathcal{ B}_{ Z }^\delta (x))}{\sigma^2_V(x/2)}\ge -1
\]
and
\[
\limsup_{x\to\infty}\frac{\ln (\mathcal{ B}_{ Z }^\delta (x))}{\sigma^2_V(x/2)}\le -\frac{3-2\sqrt{2}}{2}.
\]
\end{sat}
\BRM 
\label{rem.log}
\begin{enumerate}[(i)]
	\item If we suppose additionally that $\sigma_V^2$ is regularly varying at $\infty$ with parameter $\alpha>0$, then
 it follows   from \Cref{prop.bound} that
 \[
- \frac{1}{2^{\alpha}}\le\liminf_{x\to\infty}\frac{\ln (\mathcal{ B}_{ Z }^\delta (x))}{\sigma^2_V(x)}\le
\limsup_{x\to\infty}\frac{\ln (\mathcal{ B}_{ Z }^\delta (x))}{\sigma^2_V(x)}\le -\frac{3-2\sqrt{2}}{2^{\alpha+1}}.
\]
%Using that for $V(t)=\sqrt{2}B_{1/2}(t)$, by (\ref{b.bm}),
%$$
%\lim_{x\to\infty}\frac{\ln (\mathcal{ B}_{ Z }^\delta (x))}{\sigma^2_V(x)}=-\frac{1}{8}.
%$$
%and for $V(t)=\sqrt{2}B_{1}(t)$, by (\ref{b.dbm}),
%$$
%\lim_{x\to\infty}\frac{\ln (\mathcal{ B}_{ Z }^\delta (x))}{\sigma^2_V(x)}=-\frac{1}{8},
%$$
%we conclude that none of the obtained in Proposition \ref{prop.bound} bound is exact.
%\\
\item If follows from the proof of  \Cref{prop.bound} that under {\bf A1-A2}
 \[
- \frac{1}{2}\le\liminf_{x\to\infty}\frac{\ln (\pk{\epsilon_\delta(Y)>x})}{\sigma^2_V(x/2)}\le
\limsup_{x\to\infty}\frac{\ln (\pk{\epsilon_\delta(Y)>x})}{\sigma^2_V(x/2)}\le -\frac{3-2\sqrt{2}}{2}.
\]
\end{enumerate}
\ERM

\begin{example}
Let  $V(t)=\sqrt{2}B_H(t)$, with $H\le 1$, i.e., $\sigma^2_V(t)=2t^{2H}$.
Then $\E{\epsilon_0(Y)}=\frac{4^{1/(2H)+0.5}}{\sqrt{\pi}\Gamma(1/(2H)+0.5)}$, see \cite{KSojourn}. For $\delta>0$ we use (\ref{expepy})
to compute  $\E{\epsilon_\delta(Y)}$, see Tab. \ref{Tab:1}.
The graph of $\E{\epsilon_\delta(Y)}$ as a function of $\delta$ and the upper bound (\ref{upBMar})
with $p=1$ for Berman constants as a function of $x\in[1,10]$ are presented on Fig. \ref{Pic:1}.
\begin{table}[h!!!]
\centering
\resizebox{\textwidth}{!}{%
\begin{tabular}{l|l|l|l|l|l|l|l|l|l|l}
\hline
$\E{\epsilon_\delta(Y)}\,\,\backslash\,\, H $  & $0.1$    & $0.2$     & $0.3$ & $0.4$     &  $0.5$      &  $0.6$      &  $0.7$      & $0.8$     & $0.9$ & $1$      \\
\hline
$\delta=0$   & $60480$&   $ 72.216267$&    $12.309822$&   $  5.866446$&     $4$& $3.198992$&     $2.777685$&     $2.527405$&    $ 2.366354$&     $2.256758$\\
\hline
$\delta=1$   & $48824.040913$   & $72.594979$ &   $12.632020$  &   $6.140195$     & $4.232120$ & $3.395236$  &   $2.942920$ &  $2.665777$&    $2.481422$ &     $2.351603$\\
\hline
$\delta=5$     &   $57667.986631$ &$   74.736598 $&   $14.803952$&     $8.344951$&    $ 6.474827$&
 $5.685059$ &    $ 5.295399$&     $5.104008$&     $5.026130$&    $ 5.004070$\\
\hline
$\delta=10$  &  $59291.12614$&    $77.87128$  &  $18.22637$ &   $12.07630$&    $10.54057$& $10.09794$ &$    10.00788$ &    $10.00016$ &    $10$&    $10$\\
\hline
\end{tabular}%
}

\vspace{3mm}
\caption{Values of $\E{\epsilon_\delta(Y)}$ for $\delta=\{0,1,5,10\}$ and $V(t)=\sqrt{2}B_H(t)$.}
\label{Tab:1}
\end{table}

\begin{figure}[h!]
    \includegraphics[width=.22\textwidth, height=0.2\textheight]{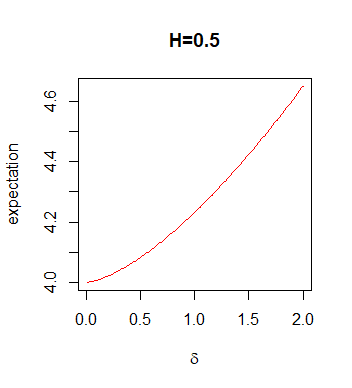}\hfill
    \includegraphics[width=.22\textwidth,height=0.2\textheight]{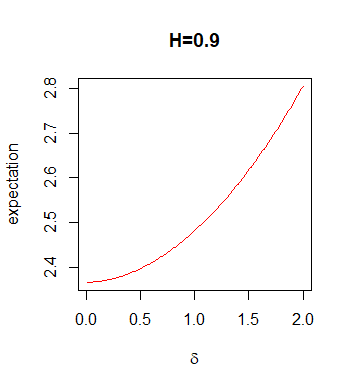}\hfill
    \includegraphics[width=.22\textwidth,height=0.2\textheight]{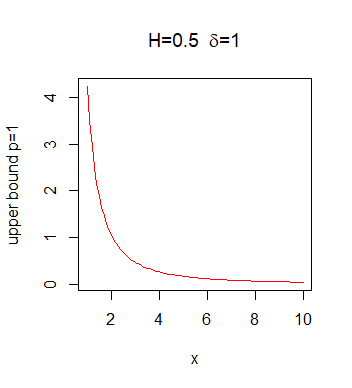} \hfill
\includegraphics[width=.22\textwidth,height=0.2\textheight]{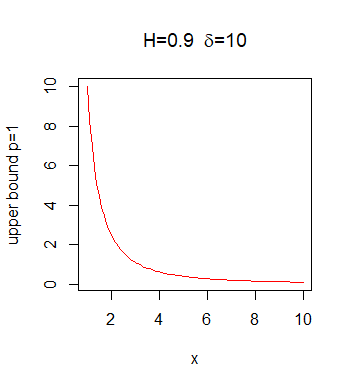}
    \caption{$\E{\epsilon_\delta(Y)}$ as a function of $\delta\in [0,2]$ and $H=\{0.5, 0.9\}$  and the upper bound (\ref{upBMar}) with $p=1$ for Berman constants as a function of $x\in[1,10]$ for $H=0.5$, $\delta=1$ and $H=0.9$, $\delta=10$  where $V(t)=\sqrt{2}B_H(t)$.}
\label{Pic:1}
\end{figure}

We simulated Berman constant $\mathcal{ B}_{ Z } (x)$ using estimator (\ref{repBE}) for different $x$ and $H$ see Tab. \ref{Tabsim2}. In our simulation we generated $N=20000$ trajectories by means of Davies-Harte algorithm on the interval $[-64, 64]$ with the step $e=1/2^9=0.001953125$. Since the sample paths of fractional Brownian motion are very torn by the negative correlation of increments for $H<0.5$ we cannot trust the simulation for $H$ close to 0 and we estimated Berman constant for $H\geq 0.4$ (see the half width of 95\% confidence interval in Tab. \ref{Tabsim2}). Let us note that the estimator (\ref{repBE}) for $x=0$ is different from the estimator of Pickands constant in \cite{DiekerY}. Compare our simulation for $x=0$ with the results of \cite{DiekerY} for Pickands constant.

\begin{table}[h!!!]
\centering
\resizebox{\textwidth}{!}{
\begin{tabular}{l|l|l|l|l|l|l|l}
\hline
$x\backslash H $  & $0.4$    & $0.5$     & $0.6$ & $0.7$     &  $0.8$ &$0.9$ & $0.999$     \\
\hline
$0$   & $1.016664$ &   $1.038244$&    $0.9336774$&   $0.8059135$&     $0.7146539$ & $ 0.6482967$ & $0.5644909$\\
        & $\pm 0.08504517$&   $\pm  0.08105526$&    $\pm  0.08679748$&   $\pm 0.04206874$&     $\pm 0.01835475$ & $\pm 0.01264331$ & $\pm  0.00675525$\\
        \hline
$0.05$   & $0.7319381$ &   $0.7662378$&    $0.77076$&   $0.7419288$&     $ 0.6972516$ & $  0.641973$ & $ 0.5660913$\\
        & $\pm 0.02164056$&   $\pm  0.02001129$&    $\pm   0.0188357$&   $\pm 0.01654269 $&     $\pm 0.01375789 $ & $\pm 0.01041449$ & $\pm  0.006183973 $\\
        \hline
$0.1$   & $ 0.6285574$ &   $ 0.6975815$&    $ 0.7083258$&   $ 0.6965521$&     $ 0.670074$ & $0.6275897 $ & $0.5638993$\\
        & $\pm  0.01509329$&   $\pm   0.01508846$&    $\pm   0.01388002$&   $\pm  0.01262112 $&     $\pm 0.01129947$ & $\pm 0.008966844 $ & $\pm  0.005561737$\\
        \hline
$0.2$   & $ 0.5210044$ &   $ 0.5855574$&    $ 0.6169511$&   $ 0.6341507$&     $ 0.6143437$ & $ 0.5951561 $ & $ 0.5601536$\\
        & $\pm 0.01018901$&   $\pm   0.01005615$&    $\pm 0.009730012 $&   $\pm 0.009250254 $&     $\pm 0.008171438$ & $\pm  0.007041064$ & $\pm   0.004994959$\\
\hline
$0.4$   & $0.3996749$ &   $ 0.4631784$&    $0.5050474$&   $0.5339913$&     $0.5503154$ & $0.5457671 $ & $ 0.5397937$\\
        & $\pm 0.00644895$&   $\pm  0.006531524$&    $\pm  0.006421914$&   $\pm  0.006178146$&     $\pm 0.005904937$ & $\pm  0.005230998$ & $\pm 0.004207075 $\\
\hline
$0.5 $  & $ 0.3523235$    & $ 0.4228112$     & $ 0.4642559$ & $ 0.4931276 $     &   $0.5079599 $ &$0.5182281$ & $ 0.5306746$     \\
  & $\pm 0.005388992$    & $\pm 0.005689292$     & $\pm 0.005635201$ & $\pm  0.005447861$     &  $\pm 0.005075071$ &$\pm 0.004643272$ & $\pm 0.003954409$     \\
\hline
$1$   & $ 0.2390332$&   $ 0.2833468$&    $ 0.3098104$&  $0.3418117$ &     $ 0.3673532$ & $0.3978701$ & $ 0.4387592$\\
   & $\pm 0.003261513$    & $\pm 0.00342165$     & $ \pm 0.003471659$ & $\pm 0.003439146$     &  $\pm 0.003371701$ &$\pm 0.003232641$ & $\pm 0.002982276$     \\
\hline
$2 $  & $0.1373955$    & $0.1493499 $     & $ 0.1606709$ & $  0.1731576$     &  $ 0.1824689$ &$0.1975992$ & $0.2061034$     \\
   & $\pm 0.001904439$&   $ \pm 0.002068819$&    $ \pm 0.002168721$& $ \pm 0.002275878$  &  $\pm 0.002382227$   & $\pm  0.002493246$ &  $\pm 0.00260955$\\
\hline
$5$   & $ 0.04407294$&   $  0.0366231$&    $0.02843742 $&  $ 0.02065171$ &   $ 0.01204482 $  & $ 0.004497369$ & $  0.001105259$\\
   & $ \pm 0.0008552848$&   $\pm  0.0008681481$&    $\pm  0.0008307949$& $\pm 0.0007567421$ &  $\pm 0.0006104972 $   & $\pm 0.0003907966$ & $  \pm 0.0001991857$\\
\hline
$6$   & $0.03259252 $&   $  0.02488601$&    $0.01711836$&  $  0.01001073$ &   $ 0.003932873$  & $ 0.0006936917$ & \rm{6.375708e-05}\\
   & $ \pm 0.0007073306$&   $\pm 0.0006875452 $&    $\pm  0.000618167$& $\pm 0.0005032682$ &  $\pm  0.0003303923$   & $\pm 0.0001440762$ & $  
\pm 4.42043\!\!$ \rm{e-05}\\
\hline
\end{tabular}
}
\vspace{3mm}
\caption{Estimation of $\mathcal{ B}_{ Z } (x)$ for fractional Brownian motion $B_H$ and the half width of 95\% confidence interval.}
\label{Tabsim2}
\end{table}

\end{example}

\COM{
\zm{I completely do not understand the simulation for fBm. They are decreasing as a function of $H$ for fixed $x=0$ and $x=5$ and they are increasing as a  function $H$ for $x=0.5, 1, 2$. It is not in the contradiction with the theoretical results but it is very strange!!! The simulation are for 5000 trajectories, $T=64$, step $e=1/2^9=0.001953125$. If I decrease $T$ by 2 or decrease $e$ by 2 my supercomputer starts to compute very long (I must interrupt it because my laptop is heating). Later we remove the results for $N=5000$.}
}

\begin{example}
Let  $X(t)$, $t\inr$ be a stationary Ornstein-Uhlenbeck process, i.e., a centered Gaussian process
with zero mean and covariance $\E{X(t) X(s)}=\exp(-|t-s|),s,t\in \R$. Then the random process
$$
V(t)=
\left\{\begin{array}{ll}
\sqrt{2}\int_0^t X(s)ds &\mbox{if $t\geq 0$}\\
&\\
-\sqrt{2}\int_t^0 X(s)ds &\mbox{if $t<0$}
\end{array}
\right.
$$
is Gaussian with  stationary increments and variance
$\sigma^2_V(t)=4(|t|+e^{-|t|}-1$).
Using (\ref{repBE}) we simulated the Berman constant for $\delta=0$ and different $x$, see Tab.\ \ref{Tabsim} and for $x=0$ and $\delta=\{0, 0.1, 0.2, 0.5, 1,2,5, 10\}$, see Tab.\ \ref{Tabsimdelta}. We generated $N=20000$ trajectories with the step $e=10^{-5}$ on the interval $[-15, 15]$.
In Fig. \ref{Pic:2}
we graphed $\mathcal{ B}_{ Z }(x)$ and $\frac{\ln (\mathcal{ B}_{ Z }(x))}{\sigma_V^2(x/2)}$ as function of $x$ and we get that this ratio is asymptotically around $-0.4$. Note that according to \Cref{rem.log} it should be between $-0.5$ and $-0.04289322$.

Using (\ref{expepy}) we computed $\E{\epsilon_\delta(Y)}$ and
$1/\E{\epsilon_\delta(Y)}$
for $\delta=\{0,0.1,0.2, 0.5,1,2,5,10\}$, see Tab. \ref{Tab:2}.
%Let us recall that $1/\E{\epsilon_\delta(Y)}$ is a lower bound for $\mathcal{ B}^{\delta}_{ Z }(0) $.
The graph of the lower bound of
$\mathcal{ B}^{\delta}_{ Z }(0) $ for the integrated Ornstein-Uhlenbeck process that is $1/\E{\epsilon_\delta(Y)}$ as a
function of $\delta\in [0,10]$ is given in Fig. \ref{Pic:2}.  The value of
$\mathcal{B}_{ Z }(0) $
constant for the integrated Ornstein-Uhlenbeck
process with the same parameters as here was simulated in \cite{demiro2003simulation} resulting in the value $0.528$.

\begin{table}[h!!!]
\centering
\resizebox{\textwidth}{!}{
\begin{tabular}{l|l|l|l|l|l|l|l}
\hline
$x $  & $0$    & $0.5$     & $1$ & $1.5$     &  $2$ &$2.5$ & $3$     \\
\hline
$\mathcal{ B}_{ Z } (x)$   & $0.5267956$&   $0.452556$&    $0.3482289$&   $ 0.2621687$&     $ 0.1900299$ & $0.1376086$ & $0.09881259$\\
   & $\pm 0.01817717$&   $\pm 0.004676632$&    $\pm 0.003180162$&   $\pm 0.002588018$&     $\pm 0.002216284$ & $\pm 0.001910763$ & $\pm 0.00163841$\\
\hline
$x $  & $4$    & $5$     & $6$ & $7$     &  $8$ &$9$ & $10$     \\
\hline
$\mathcal{ B}_{ Z } (x)$   & $ 0.05088893$&   $0.02715927$&    $0.01433577$&   $0.007437053$&     $ 0.003796336$ & $0.001998398$ & $0.001205136$\\
   & $\pm 0.00116684$&   $\pm 0.0008278098$&    $ \pm 0.0005788133$&   $\pm 0.0003983809$&     $\pm 0.0002730899$ & $\pm 0.0001906838$ & $\pm 0.0001414664$\\
\hline
$x $  & $11$    & $12$     & $13$ & $14$     &  $15$ &$16$ & $17$     \\
\hline
$\mathcal{ B}_{ Z } (x)$   & $ 0.000631948$&  $0.0003812784$ &    $0.0001845301$&   \rm{0.00010499}&     {\rm 9.130422e-05} & {\rm 2.426165e-05} & {\rm 2.103512e-05}\\
   & $\pm 9.837308\!\!$ \rm{e-05}&   $\pm 7.355149\!\!$ \rm{e-05}&    $\pm 4.89203\!\!$ \rm{e-05}&   $\pm 3.593126\!\!$ \rm{e-05}  &   $\pm 3.276786\!\!$ \rm{e-05}   & $\pm 1.594309\!\!$ \rm{e-05} & $\pm  1.463212e\!\!$ \rm{e-05}\\
\hline
\end{tabular}
}
\vspace{3mm}
\caption{Estimation of $\mathcal{ B}_{ Z } (x)$ for integrated Ornstein-Uhlenbeck process and the half width of 95\% confidence interval..}
\label{Tabsim}
\end{table}

\begin{table}[h!!!]
\centering
\resizebox{\textwidth}{!}{
\begin{tabular}{l|l|l|l|l|l|l|l|l}
\hline
$\delta $ &$0$ & $0.1$    & $0.2$     & $0.5$ & $1$     &  $2$ &$5$ & $10$     \\
\hline
$\mathcal{ B}^{\delta}_{ Z }(0) $ &$0.5267956$  & $0.5131973$&   $ 0.5126575$&    $ 0.4934096$&   $ 0.4484668$&     $0.3843544$ & $0.1908583 $ & $ 0.09984$\\

  &$\pm 0.01817717$  & $\pm 0.007850686$&   $\pm 0.007194036$&    $\pm 0.005746635$&   $\pm 0.00402201$&     $ \pm 0.001995524$ & $\pm 0.0004065713$ & $\pm 3.913749\!\!$ \rm{e-05}\\
\hline
\end{tabular}
}
\vspace{3mm}
\caption{Estimation of $\mathcal{ B}^{\delta}_{ Z }(0) $ for integrated Ornstein-Uhlenbeck process and the half width of 95\% confidence interval.}
\label{Tabsimdelta}
\end{table}

\begin{table}[h!!!]
\centering
\resizebox{\textwidth}{!}{
\begin{tabular}{l|l|l|l|l|l|l|l|l}
\hline
$\delta $  & $0$    & $0.1$     & $0.2$ & $0.5$     &  $1$&$2$ & $5$ & $10$       \\
\hline
$\E{\epsilon_\delta(Y)}$   & $3.234658$&   $3.245584$&    $3.248405$&   $3.268183$&     $3.339158$ &$3.626068$ & $5.482154$  & $10.05426$\\
\hline
$1/\E{\epsilon_\delta(Y)}$   & $0.3091517$&   $0.308111$&    $0.3078434$&   $0.3059804$&     $0.2994767$ &$0.2757808$ & $0.1824101$  & $0.09946035$
\\
\hline
\end{tabular}
}
\vspace{3mm}
\caption{Estimation of $\E{\epsilon_\delta(Y)}$ and $1/\E{\epsilon_\delta(Y)}$ for $\delta=\{0,0.1,0.2, 0.5,1,2,5,10\}$.}
\label{Tab:2}
\end{table}
%The graph of $\E{\epsilon_\delta(Y)}$ as function of $\delta\in [0,10]$ and

\begin{figure}[h!]
\includegraphics[width=.3\textwidth,height=0.25\textheight]{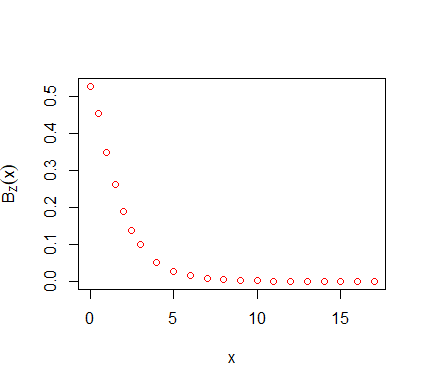}\hfill
\includegraphics[width=.3\textwidth,height=0.25\textheight]{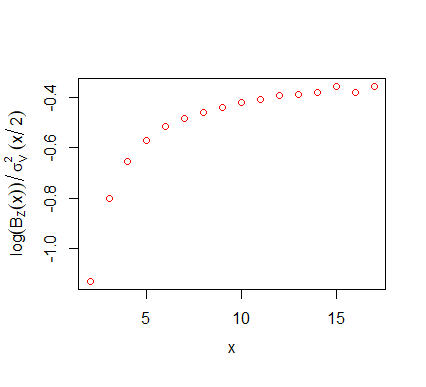}\hfill
\includegraphics[width=.3\textwidth, height=0.25\textheight]{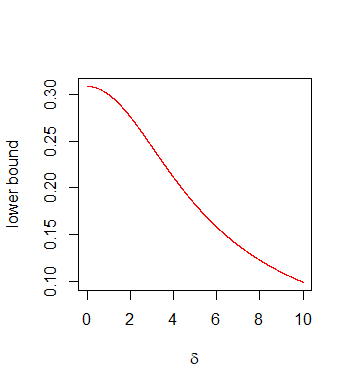}
    \caption{The graphs of $\mathcal{ B}_{ Z }(x)$ and
%$\E{\epsilon_\delta(Y)}$ as function of $\delta$
$\frac{\ln (\mathcal{ B}_{ Z }(x))}{\sigma_V^2(x/2)}$ as function of $x$
and the lower bound of Pickands constant as a function of $\delta$ for integrated Ornstein-Uhlenbeck process.}
\label{Pic:2}
\end{figure}
\end{example}

\section{Further Results and Proofs}
 Let in the following $X(t),t\inr^d$ be a max-stable stationary rf with \cadlag sample paths
and spectral rf $Z$ as in \netheo{erdha} and define $\Theta$ as in \eqref{thet}.
Define $Y(t)=R \Theta(t),t\inr^d$  with $R$ an $1$-Pareto rv \zm{(with survival function
	$1/x, x\ge 1 $)} independent of $\Theta$ and set  $M_{Y,\delta }= \sup_{t\in \delta \Z^d} Y(t)$. Note in passing that
$$\pk{ M_{Y, \delta } > 1}=1$$
 since $M_{Y, \delta } \ge  R \Theta(0)> 1$ almost surely by the assumption on $R$ and  by the definition
$\pk{\Theta(0)=1}=1$.

 Recall that $S_\delta=S_\delta(Z)= \int_{\delta \Z^d}  Z (t) \lambda_\delta(dt)$. In view of \eqref{sid} we have that $\pk{ S_0>0}=1$.
In the following, for any fixed $\delta \ge 0$ (but not simultaneously  for two different $\delta$'s) we shall assume that $S_\delta >0$ almost surely, i.e., $Z$ is such that $\pk{ \sup_{t\in \delta \Z^d} Z(t)> 0}=1$.
Such a choice of $Z$ is possible   in view of \cite[Lem 7.3]{HBernulli}. \\
A functional $F:D \to [0,\IF]$ is  said to be shift-invariant if $F(f(\cdot -h) )= F(f(\cdot))$
for all   $h\in \R^d$.

We state first two lemmas and proceed with the postponed proofs.

\BEL \label{L1} If $\pk{ S_0 < \IF}=1$,  then $\pk{ S_\delta< \IF}=1, \delta >0$  and for all  $x>0$
\bqn{ \label{pus}
	\pk{\epsilon_\delta(Y/x)< \IF}=1, \quad  \forall x>0, \forall \delta \ge 0.
}
Moreover $\pk{ M_{Y, \delta}< \IF}=1$.
\EEL
\prooflem{L1}
In view of  \eqref{sid} $S_0>0$ almost surely. The assumption that $S_0< \IF$ almost surely is in view of \cite[Thm 3]{dom2016}
equivalent with  $Z(t) \to 0$ almost surely as $\norm{t}\to \IF$, with $\norm{\cdot}$ some norm on $\R^d$. Hence $S_\delta < \IF$ almost surely follows from  \cite[Thm 3]{dom2016}.  By the definition of $\Theta$ and \cE{the fact that $\pk{ S_0(Z) \in (0,\IF)}=1$} we have
\bqn{ \label{NT}
	0
%	&=&    \E*{\ind{ \liminf_{ \abs{t} \to \IF}   Z (t)>0}} =    \E*{ \frac{S_0(Z)}{S_0(Z)} \ind{ \liminf_{ \abs{t} \to \IF}    Z (t)>0}  }\notag \\
%	&=&  \intT \E*{ \frac{  Z (h)}{S_0(Z)} \ind{ \liminf_{ \abs{t} \to \IF}   Z (t) >0} } \lambda(dh) \notag\\
	&=&  \cE{  \E*{ {  Z (0)}  \ind{ \limsup_{ \norm{t} \to \IF}  Z(t) >0} }}
	=  \cE{  \E*{  Z (0)\ind{ \limsup_{ \norm{t} \to \IF}  Z(t) /Z(0)>0} }} \notag\\
%	&=&\kdd{??}  \intT \E*{ \frac{  Z (0)}{S_0(Z/Z(0))} \ind{ \liminf_{ \abs{t} \to \IF}   \frac{ Z (t-h)}{ Z (0)} >0} } \lambda(dh) \notag\\
	&=&   \E*{  \ind{ \limsup_{ \norm{t} \to \IF}  \Theta(t) >0} }
	%\notag \\
%	&=&  \intT \E*{ \frac{ 1 }{S_0(\Theta)} \ind{ \liminf_{ \abs{t} \to \IF}   \Theta^\alpha(t)>0} } \lambda(dh)
}
implying that $\pk{ \limit{\norm{t}}  \Theta(t) =0}=1$.  Consequently, $\pk{ \limit{\norm{t}}  Y(t) =0}=1$ and hence
 the claim follows.

\QED

Below we interpret $\IF \cdot 0$ and $0/0$ as 0. The next result is a minor extension of \cite[Lem 2.7]{PH2020}.
\BEL\label{ava} If $\pk{S_0<\IF}=1$, then
for all  measurable shift invariant functional $F$ and all $ \delta,x$ non-negative
\bqn{ x  \E*{ \frac{F(Y/x)}{ \epsilon_\delta(Y)}\ind{ M_{Y,\delta}>  {\max}(x,1)} } &=&
\E*{ \frac{F(Y)}{ \epsilon_\delta(Y)}\ind{M_{Y,\delta}>  {\max}(1/x,1) }}.
}
\EEL
\prooflem{ava}   For all measurable functional $F: D \to [0,\IF]$
and all $ x> 0$
\bqn{  x \E{  F(Y) \mathbb{I}( { Y(h)} >x )   }&=&
\E{ F(x B^h Y) \mathbb{I}(x { Y(-h)} >1 )}
	\label{tYY}}
is valid for all $h \in \TT$ with $B^h Y(t)=  Y(t-h), h,t\inr^d$. \cE{Note in passing that
	$B^h Y$ can be substituted by $Y$ in the right-hand side of \eqref{tYY} if $F$ is shift-invariant.} The identity
\eqref{tYY} is shown in  \cite{MartinE}. For the discrete setup it is shown initially in \cite{Hrovje,BP} and for case $d=1$ in \cite{PH2020}.\\
Next, if $x \in (0,1]$, since $Y(0)=R> 1$ almost surely and by the assumption on the sample paths we have that
$\pk{\epsilon_\delta (Y/x)>0}=1,$  recall $\pk{\Theta(0)=1}=1$.   By \Cref{L1}
$\pk{ M_{Y, \delta} \in (1, \IF)}=1$, hence for all $x> 1$ we have further that $ M_{Y, \delta } > x$ implies $\epsilon_\delta (Y/x)>0$. Consequently, in view of \eqref{pus}
$\epsilon_\delta (Y/x)/\epsilon_\delta (Y/x)$
is well defined on the event $M_{Y, \delta}  > x, x> 1$ and also it is well-defined for any $x\in (0,1]$.

Recall that $\lambda_\delta(dt)$ is the Lebesgue measure on $\R^d$ if $\delta=0$ and the counting measure multiplied by $\delta^d$ on $\delta \Z^d$ if $\delta >0$.
{Let us remark that \zm{for any shift-invariant functional $F$}, the functional
$$
F^*(Y)=\frac{F(Y/x)\ind{ M_{Y,\delta}> \max(x,1) }}{\epsilon_\delta(Y)\epsilon_\delta(Y/x)}
$$
is shift-invariant for all $h\in\R^d$ if $\delta=0$ and any shift $h\in\delta\Z^d$ if $\delta>0$.}
Thus applying the Fubini-Tonelli theorem twice and \eqref{tYY} with functional $F^*$ we obtain  for
all $\delta \ge 0,x>0$
\bqny{ \lefteqn{
		x\E*{ \frac{F(Y/x)}{\epsilon_\delta(Y)} \ind{M_{Y,\delta}> \max(x,1) }}} \\
	&=&
	x\int_{ \delta \Z^d}
	\E*{\frac{F(Y/x)\ind{ M_{Y,\delta}> \max(x,1) }}{\epsilon_\delta(Y)\epsilon_\delta (Y/x) } \ind{Y(h)> x}}  \lambda_\delta(dh)\\
	&=& 	\int_{ \delta \Z^d}
	\E*{\frac{F(Y)\ind{ M_{Y,\delta}> \max(1/x,1) }}{\epsilon_\delta(xY)\epsilon_\delta (Y) } \ind{xY(-h)> 1}}  \lambda_\delta(dh)\\
&=& 	\E*{ \frac{F(Y)\ind{ M_{Y,\delta}> \max(1/x,1) }}{\epsilon_\delta(xY)\epsilon_\delta (Y) } \int_{ \delta \Z^d}\ind{xY(h)> 1} \lambda_\delta(dh)} \\
&=&
\E*{ \frac{F(Y)}{ \epsilon_\delta(Y)}\ind{ M_{Y,\delta}> \max(1/x,1) }},
}
hence the proof follows. \QED

\prooftheo{erdha}   Let $\delta \ge 0$ be fixed and consider for simplicity $d=1$. By the assumption we have
$\E{ \sup_{t\in [0,T]}  Z (t)}< \IF$ for all $T>0$. Since we assume that $\pk{\sup_{t \in \R} Z(t)>0}=1$, then $\pk{S_0>0}=1$. Using the assumption we have $\pk{ S_\eta< \IF}>0$  for all $\eta\ge 0$ and thus by  \eqref{rinashero}    we obtain
%\zm{WE DEFINED $S_0$ WITHOUT ANY $\alpha$ ABOVE}
\bqny{
	\IF   & > &   \E*{ \sup_{t\in [0,\bE{2}]  } Z (t) \frac{S_0}{S_0}}	= \int_{\R}\E*{ Z (h) \sup_{t\in [0,2]  }  Z (t) /S_0} \lambda(dh)\\
	& =&
\int_{\R}   \E*{ Z (0) \sup_{t\in [-h,2-h]  }  Z (t) /S_0} \lambda(dh)\\
& =&
\zm{\sum_{i\in \Z}  \int_i^{i+1}\E*{ Z (0) \sup_{t\in [-h, 2-h]  }  Z(t) /S_0}  \lambda(dh)}\\
	& {\geq }&
\sum_{i\in \Z}  \E*{ Z(0) \sup_{t\in [-i,{1}-i]  }  Z(t) /S_0} \\
	&\ge &\E*{ Z(0) \sup_{t\in \R  }  Z(t) /S_0 }. }
Since also for any $M>0$ and $\eta>0$
\bqny{
0 <  \E*{ \sup_{t\in [0,{M}]  } Z (t) \frac{S_\eta}{S_\eta}} < \IF
}	
we conclude as above that for all  $\eta\ge 0$
\bqn{\label{meteo}
	\E*{ Z(0)  \sup_{t\in \R} Z (t) /S_\eta  }< \IF, \quad
		\int_{\eta \Z} \E*{  Z (0) \sup_ {h \le t\le h +1} Z (t) /S_\eta }\lambda_\eta(dh)  < \IF .
}
Next, for any $x\ge 0$ and $\eta \ge 0$
\bqny{ \lefteqn{
		 T^{-1}  \int_0^\IF \pk{    \int_{[0,T] \cap \eta \Z} \mathbb{I}( Z  (t) >   s) \lambda_\eta(dt) > x ,  S_\eta=\IF} ds}\\
	&\le &	   T^{-1}   \int_0^\IF    \pk*{ \int_{[0,T] \cap \eta  \Z} \mathbb{I}( Z  (t) >   s) \lambda_\eta(dt) > 0 ,  S_\eta=\IF}  ds \\
	&= &	   T^{-1}   \int_0^\IF    \pk*{ \sup_{t\in [0,T] \cap \eta \Z }  Z  (t) >   s ,  S_\eta=\IF}  ds \\
%&\le  &	  T^{-1}   \int_0^\IF    \pk*{ \sup_{t\in [0,T] \cap \eta \Z }  Z  (t) >   s  ,  S_\eta=\IF}  ds \\
&\zm{=}&   T^{-1}    \E*{  \sup_{t\in [0,T] \cap \eta \Z }  Z  (t) , S_\eta=\IF }\\
&\to & 0, \quad T \to \IF,
}	
%\zm{DO WE NEED THE LAST MAX THAT IS THE FIRST TERM IS NOT ENOUGH?}
where the last claim follows from \cite[Cor 2.1]{ZKE}.
We shall assume that $S_\delta >0$ almost surely (this is possible as mentioned at the beginning of this section). For notational simplicity we consider next $\delta=0$.
For any $M>0, T> 2M$ by the Fubini-Tonelli theorem \zm{and} \eqref{rinashero}
\bqny{
\lefteqn{\int_0^\IF \pk*{    \int_{[0,T] } \mathbb{I}( Z  (t) >   s) \lambda(dt) > x  } ds}   \\
	&=&    \int_{\R} \E*{  \frac{ Z  (h)}{S_{\zm{0}}} \int_0^\IF \ind{ \int_{[0,T] } \mathbb{I}( Z  (t) >   s) \lambda(dt) > x    } ds }  dh  \\
	&=&    \int_{\R} \E*{  \frac{ Z  (0)}{\cE{S_{\zm{0}}}} \int_0^\IF\ind{  \int_{[0,T]} \mathbb{I}( Z  (t-h) >   s) \lambda(dt) > x }ds} dh\\
	&=&    \int_{\R} \int_0^\IF \E*{  \frac{ Z  (0)}{S_{\zm{0}}}  \ind{ \int_{h}^{ T+h} \mathbb{I}( Z  (t) >   s) \lambda(dt) > x }} dsdh\\	
&=&    \int_{-M-T}^{-M}\int_0^\IF \E*{  \frac{ Z  (0)}{S_{\zm{0}}}  \ind{ \int_h^{ T+h} \mathbb{I}( Z  (t) >   s) \lambda(dt) > x }} dsdh\\
&& + 	   \int_{h<-M-T \,\text{or}\, h>-M}\int_0^\IF \E*{  \frac{ Z  (0)}{S_{\zm{0}}}  \ind{ \int_h^{ T+h} \mathbb{I}( Z  (t) >   s) \lambda(dt) > x }} dsdh\\
&=:&I_{M,T}+ J_{M,T}.
}	
Thus we obtain
\bqny{   \frac{I_{M,T}}{T}
	&=&  \int_{-M/T-1}^{-M/T}  \E*{ \frac{ Z (0)}{S_{\zm{0}}} \int_0^\IF  \mathbb{I} \Bigl(  \int_{Th}^{T(1+h)}  \mathbb{I}(  Z   (t) > s)  \lambda(dt) > x \Bigr)  ds } dh\\
	&=&  \int_{-1}^{0}  \E*{ \frac{ Z (0)}{S_{\zm{0}}} \int_0^\IF  \mathbb{I} \Bigl(  \int_{Th-M}^{T(1+h)-M}  \mathbb{I}(  Z   (t) > s)  \lambda(dt) > x \Bigr)  ds} dh\\
	&\to & \int_{-1}^{0}  \E*{  \frac{ Z (0)}{S_{\zm{0}}}  \int_0^\IF  \mathbb{I} \Bigl(  \int_{\delta\Z}  \mathbb{I}(  Z   (t) > s)  \lambda(dt) > x \Bigr) ds } dy, \quad T\to \IF\\
	&=&  \E*{  \frac{ Z (0)}{S_{\zm{0}}} \int_0^\IF  \mathbb{I} \Bigl(  \int_{\delta\Z}  \mathbb{I}(  Z   (t) > s ) \lambda_\delta(dt) > x \Bigr) ds }\\
	&\le &  \E*{  \frac{ Z (0)}{S_{\zm{0}}} \int_0^\IF  \mathbb{I} \Bigl(  \int_{\R}  \mathbb{I}(  Z   (t) > s ) \lambda(dt) > 0 \Bigr) ds}\\
	&= &  \E*{  \frac{ Z (0)}{S_{\zm{0}}}   \sup_{t\in \R}  Z (t)}	= \zm{\mathcal{H}_{ Z }^{\zm{0}}	}< \IF,
}
where $\zm{\mathcal{H}^{\zm{0}}_{ Z }}$ is the Pickands constants, see \cite[Prop 2.1]{ZKE} for the last formula.\\
Let us consider the second term
\bqny{\frac{J_{M,T}}{T}  &=&
 \int_{(-\infty, -1)\cup (0,\infty)}  \E*{ \frac{ Z (0)} {S_0} \int_0^\IF  \mathbb{I} \Bigl(  \int_{Th-M}^{T(1+h)-M}  \mathbb{I}(  Z   (t) > s)  \lambda(dt) > x \Bigr) ds  } dh\\
	&=& 	
	\int_{-\infty}^{-1}  \E*{ \frac{ Z (0)} {S_0} \int_0^\IF  \mathbb{I} \Bigl(  \int_{Th-M}^{T(1+h)-M}  \mathbb{I}(  Z   (t) > s)  \lambda(dt) > x \Bigr) ds } dh\\
&& + \int_0^\infty  \E*{ \frac{ Z (0)} {S_0} \int_0^\IF  \mathbb{I} \Bigl(  \int_{Th-M}^{T(1+h)-M}  \mathbb{I}(  Z   (t) > s)  \lambda(dt) > x \Bigr) ds } dh\\
	&=:& K_{M,T}+L_{M,T}\,.
}
Further,
assuming for simplicity that $T$ is a positive integer we get
\bqny{
K_{M,T}&\leq&\int_{-\infty}^{-1}  \E*{ \frac{ Z (0)} {S_0} \int_0^\IF  \mathbb{I} \Bigl(  \int_{Th-M}^{T(1+h)-M}  \mathbb{I}(  Z   (v) > s)  \lambda(dv) > 0 \Bigr) ds } dh\\
&=&\int_{-\infty}^{-1}  \E*{ \frac{ Z (0)} {S_0}  \sup_{t\in [Th-M,\, T(1+h)-M]} Z (t)} dh\\
&=&\frac{1}{T}\int_{-\infty}^{-T-M}  \E*{ \frac{ Z (0)} {S_0}  \sup_{t\in [h,\, h+T]} Z (t)} dh\\
&\leq&\frac{1}{T}\sum_{i=1}^{T-1}\int_{-\infty}^{-T-M}  \E*{ \frac{ Z (0)} {S_0} \sup_{t\in [h+i,\, h+i+1]} Z (t)} dh\\
&=&\frac{1}{T}\sum_{i=1}^{T-1}\int_{-\infty}^{-T-M+i}  \E*{ \frac{ Z (0)} {S_0} \sup_{t\in [h,\, h+1]} Z (t)} dh\\
&\leq&\frac{1}{T}\sum_{i=1}^{T-1}\int_{-\infty}^{-M}  \E*{ \frac{ Z (0)} {S_0} \sup_{t\in [h,\, h+1]} Z (t)} dh\\
&=&\int_{-\infty}^{-M}  \E*{ \frac{ Z (0)} {S_0} \sup_{t\in [h,\, h+1]} Z (t)} dh
\to 0, \quad M \to \IF,
}
where the last convergence follows from \eqref{meteo}. The same way we show that $L_{M,T}\to 0$ as $M\to\IF$ establishing the proof.\\
We prove next the second claim. In view of \cite[proof of \zm{Prop 2.1}]{ZKE}
almost surely for all $\delta,\eta \in [0,\IF)$
\bqn{ \label{lime2}
	\frac{ 1 }{S_\eta(\Theta)} =\frac{ 1 }{S_\eta(\Theta)} \frac{ S_\delta(\Theta)} {S_\delta(\Theta)}\Theta(0), \quad \{S_\eta(\Theta) <\IF \} =
	 \{S_\delta(\Theta) <\IF \}.
}	
Consequently, for any $\delta,\eta, x$ non-negative
\bqny{
\lefteqn{ B_{\delta, \eta}(x)}\\
&:=& \int_0^\IF \E*{ \frac{ Z (0)}{ S_\eta} \ind{ S_\eta< \IF}\mathbb{I} \Bigl(  \int_{\delta \Z}  \mathbb{I}(  Z   (t) >  s)  \lambda_\delta(dt) > x \Bigr) }ds\\
&=&
\int_0^\IF \E*{ \frac{1 }{ S_\eta(\Theta)} \ind{ S_\eta(\Theta)< \IF} \mathbb{I} \Bigl(  \int_{ \delta \Z}  \mathbb{I}( \Theta  (t) >  s)  \lambda_\delta(dt) > x \Bigr) }ds\\
&=&
\int_0^\IF \E*{ \frac{1 }{ S_\eta(\Theta)} \frac{S_\delta(\Theta)}{S_\delta(\Theta)} \ind{ S_\delta(\Theta)< \IF} \mathbb{I}
	\Bigl(  \int_{ \delta \Z}  \mathbb{I}( \Theta  (t) >  s)  \lambda_\delta(dt) > x \Bigr) }ds.
}
We proceed next with the case $\delta=0$, the other case follows with the same argument \zm{where} it is important that $\eta=k\delta$ for the shift transformation. Taking $\delta=0, \eta>0$ we have
\bqny{
	\lefteqn{ B_{0, \eta}(x)}\\
	&=& \int_{\R} 	\int_0^\IF \E*{ \frac{1 }{ S_\eta(\Theta)} \frac{\Theta(r)}{S_0(\Theta)} \ind{ S_0(\Theta)< \IF} \mathbb{I}
		\Bigl(  \int_{ \R}  \mathbb{I}( \Theta  (t) >  s)  \lambda(dt) > x \Bigr) }ds\lambda(dr)\\
	&=&   \sum_{v\in \eta \Z} \int_{ r+v\in [0,\eta ]} 	\int_0^\IF \E*{ \frac{ Z (0) }{ S_\eta(Z)} \frac{  Z (r)}{S_0(Z)} \ind{ S_0(Z)< \IF} \mathbb{I}
		\Bigl(  \int_{ \R}  \mathbb{I}( Z   (t) >  s)  \lambda(dt) > x \Bigr) }ds\lambda(dr)\\
	&=&   	\int_0^\IF \E*{  \sum_{v\in \eta \Z} \int_{ r +v\in [0,\eta]}\frac{ Z (0) }{ S_\eta(Z)} \frac{  Z (r)}{S_0(Z)} \ind{ S_0(Z)< \IF} \mathbb{I}
	\Bigl(  \int_{ \delta \Z}  \mathbb{I}( Z   (t) >  s)  \lambda(dt) > x \Bigr) }\lambda(dr) ds\\
	&=&   	\int_0^\IF \E*{  \sum_{v\in \eta \Z} \int_{ r\in [0,\eta]}\frac{ Z (0) }{ S_\eta(Z)} \frac{  Z (r-v)}{S_0(Z)} \ind{ S_0(Z)< \IF} \mathbb{I}
	\Bigl(  \int_{ \R}  \mathbb{I}( Z   (t) >  s)  \lambda(dt) > x \Bigr) }\lambda(dr) ds\\
	&=&   	\int_0^\IF \E*{   \int_{ r\in [0,\eta]} \frac{1}{\eta} \sum_{v\in \eta \Z} \frac{ \eta  Z (v-r) }{ S_\eta({B^r}Z)}  \lambda(dr) \frac{  Z (0)}{S_0(Z)} \ind{ S_0(Z)< \IF} \mathbb{I}
	\Bigl(  \int_{ \R}  \mathbb{I}( Z   (t) >  s)  \lambda(dt) > x \Bigr) } ds\\
	&=&   	\int_0^\IF \E*{   \frac{  Z (0)}{S_0(Z)}  \mathbb{I}
	\Bigl(  \int_{ \R}  \mathbb{I}( Z   (t) >  s)  \lambda(dt) > x \Bigr) }ds,
}
where we used \eqref{rinashero} \zm{with $h=r-v$} to obtain the second last equality above and \eqref{laps} to get the last equality,	hence the proof follows.

 \QED

\proofkorr{K1}
Given $x\ge 0$ consider  the representation \eqref{bli}
$$
\mathcal{B}_{ Z }^0(x) = \int_0^\IF \E*{ \frac{ Z (0)}{ S_0} \mathbb{I} \Bigl(  \int_{\R^d}  \mathbb{I}(  Z   (t) >  s)  \lambda(dt) > x \Bigr) }ds.
$$
By the monotonicity with respect to variable $x$ of the function
\begin{equation}\label{expxg}
\E*{ \frac{ Z (0)}{ S_{0}} \mathbb{I} \Bigl(  \int_{\R^d}  \mathbb{I}(  Z   (t) >  s)  \lambda(dt) > x \Bigr) }
\end{equation}
in order  to show the continuity of $\mathcal{B}_{ Z }^0(x)$ it suffices to prove that
\begin{equation}\label{expx}
\E*{ \frac{ Z (0)}{ S_{0}} \mathbb{I} \Bigl(  \int_{\R^d}  \mathbb{I}(  Z   (t) >  s)  \lambda(dt) = x \Bigr) }=0
\end{equation}
for almost all $s>0$.
Let us define the following measurable sets
$$
A_s=\mathbb{I} \Bigl(  \int_{\R^d}  \mathbb{I}(  Z   (t)>s)  \lambda(dt) = x \Bigr)\,.
$$
Since $Z$ has almost surely continuous trajectories we have $A_s\cap A_{s'}=\emptyset$ if $0<s<s'$ and $x>0$. Thus
there are countably many $s>0$ such that $\pk{A_s}>0$ because if there were not countably many ones we would find countably many disjoint $A_s$ such that $\sum \pk{A_s}=\IF$. Thus we get \eqref{expx} for almost all $s>0$. The continuity at $x=0$ follows from the right continuity of \eqref{expxg}. \QED

\prooflem{xhxh}  \underline{\Cref{en:A}}: In view of \eqref{bli} and  substituting $\Theta(t)=Q_\delta(t)/ S_\delta(\Theta)$ to \eqref{eqTE} we get
$$
\mathcal{B}_{ Z }^\delta(x)=
 \int_0^\IF \E*{ \frac{1}{S_\delta(Q_\delta)}\mathbb{I} \Bigl(  \int_{\delta \Z^d}  \mathbb{I}( Q_\delta   (t) >  s)  \lambda_\delta(dt) > x \Bigr) }ds.
$$
Since $S_\delta(Q_\delta)=1$ the claim follows.\\
\underline{\Cref{en:B}}: If   $ \pk{ S_0< \IF} > 0$ we can define
$V(t)= Z(t) \lvert { S_0 < \IF}$ and set (recall $S_0=S_0(Z), \E{  Z (0)}=1$)
$$b=  \E{  Z (0)\ind{ S_0 < \IF}}= \pk{S_0(\Theta)< \IF  }>0.$$
For this choice of $b$ {by \eqref{rinashero}} we have
$$\E{V (t)} = \frac{\E{ Z(t) \ind{ S_0 < \IF} }}{ \pk{S_0(\Theta)< \IF  }} =
\frac{\E{ Z(0) \ind{ S_0 < \IF} }}{ \pk{S_0(\Theta)< \IF  }}= 1$$
 for all $t\in \R$. Clearly, $\pk{\sup_{t\in \TT} V(t)>0} =1$. In view of \cite{dom2016} $V$ is the spectral rf of a  stationary max-stable rf $X_*$  with c\`{a}dl\`{a}g sample paths and moreover $S_0(V)=\intT V (t) \lambda(dt)<\IF$ almost surely.
In view of \cite[\zm{proof of Prop 2.1}]{ZKE} we have that
$$\{S_\delta(\Theta)<\IF\}= \{S_0(\Theta)<\IF\}$$
 almost surely for all $\delta>0$. Consequently,  we obtain for all  $\delta>0$
\bqny{
	\mathcal{ B}_{ Z }^ \delta(x)
	&=& \int_0^\IF \E*{ \frac{ Z (0)}{ S_\delta } \mathbb{I} \Bigl(  \int_{\delta \Z^d}  \mathbb{I}(  Z   (t) >  s)  \lambda_\delta(dt) > x \Bigr) \ind{ S_\delta< \IF }}ds\\
	&=& \int_0^\IF \E*{ \frac{ 1}{ S_\delta(\Theta) } \mathbb{I} \Bigl(  \int_{\delta \Z^d}  \mathbb{I}(  \Theta   (t) >  s)  \lambda_\delta(dt) > x \Bigr) \ind{ S_\delta(\Theta)< \IF, S_0(\Theta) < \IF} }ds\\
	&=& \int_0^\IF \E*{ \frac{ Z (0)}{ S_\delta(Z) } \mathbb{I} \Bigl(  \int_{\delta \Z^d}  \mathbb{I}(  Z   (t) >  s)  \lambda_\delta(dt) > x \Bigr) \ind{ S_0(Z)< \IF} } ds\\
	&=& b \int_0^\IF \E*{ \frac{ V(0)}{ S_\delta(V) } \mathbb{I} \Bigl(  \int_{\delta \Z^d}  \mathbb{I}(  V (t) >  s)  \lambda_\delta(dt) > x \Bigr) }ds\\
	&=& 	b\mathcal{ B}_{ V}^ \delta(x) < \IF
	\label{zhgalla}
}
establishing the proof.
\QED

\prooftheo{phl} Assume first that $\pk{S_0< \IF}=1$. In view of \eqref{pus} we have that $\epsilon_\delta < \IF$ almost surely, hence as in  \cite{booksoulier, PH2020} where $d=1$ is considered  it follows that  \eqref{Q} holds with
$$
\cE{\bar Q_\delta (t)}= c \frac{ Y(t)}{ \epsilon_\delta(Y)    M_{Y,\delta }} , \quad t\in \R,
$$
\cE{with $c=1$ if $\delta=0$ and $c=\delta^d $ otherwise}. Set below $Q_\delta= \bar Q/c$
and for simplicity omit the subscript below writing  simply $M_Y$ instead of $M_{Y,\delta}$.
Since $Y(t)/M_Y \le 1$  almost surely for all $t\in \delta \Z^d$ and $\pk{ M_Y\in (1,\IF)}=1$,  in view of \nelem{xhxh} we have using further the Fubini-Tonelli theorem and \nelem{ava}
\bqny{
	 	\mathcal{ B}_{ Z }^\delta (x)  &=& \int_0^\IF \E*{  \mathbb{I} \Bigl(  \int_{\delta \Z^d}  \mathbb{I}(  {Q} _\delta   (v) >  s)  \lambda_\delta(dv) > x \Bigr) }ds  \\
%	\int_0^\IF  \E*{\frac{1}{\epsilon_\delta (Y) }       \mathbb{I} \Bigl(  \int_{\delta \Z^d}  \mathbb{I}( Y(v) /M_Y   > s )  \lambda_\delta(dv) > x \Bigr)  } ds\\
%	&=& 	\int_0^\IF  \E*{\frac{1}{\epsilon_\delta (Y)  }  \frac{1}{M_Y  }      \mathbb{I} \Bigl(  \int_{\delta \Z^d}  \mathbb{I}( Y  (v) >  s)  \lambda_\delta(dv) > x \Bigr)  } ds\\
	&=& 	\int_0^\IF  \E*{\frac{1}{\epsilon_\delta (Y)  }  \frac{1}{M_Y  }  \ind{ \epsilon_\delta (Y/s) > x}} ds\\
	&=& 	\int_0^\IF  \E*{\frac{1}{\epsilon_\delta (Y)  }  \ind{M_Y > s} \frac{1}{M_Y  }  \ind{ \epsilon_\delta (Y/s) > x}} ds\\
	&=:& 	\int_0^\IF  \E*{\frac{1}{ {s}\epsilon_\delta (Y)  }  \ind{M_Y > s} F(Y/s)  } ds\\
&=& 	\int_0^\IF  \frac{1}{s^2}  \E*{  \frac{ \ind{\epsilon_\delta (Y) > x} }{\epsilon_\delta (Y ) M_Y  }      \ind{M_Y > \max(1/s,1)} }  ds\\
	&=& \E*{\frac{ \ind{ \epsilon_\delta (Y) > x}}{\epsilon_\delta (Y ) M_Y  }  	      	\int_0^\IF  \frac{1}{s^2}   \ind{M_Y > \max(1/s,1)}   ds}\\
	&=& \E*{\frac{ \ind{ \epsilon_\delta (Y) > x}}{\epsilon_\delta (Y )}       }.
}	
The last equality follows from (recall $M_Y\in (1,\IF)$ almost surely)
$$\int_0^\IF  \frac{1}{s^2}   \ind{M_Y > \max(1/s,1)}   ds=
\int_1^\IF  \frac{1}{s^2}      ds+\int_0^1  \frac{1}{s^2}   \ind{M_Y > 1/s}   ds
=M_Y.$$
In view of  \eqref{pus}  for  all $x$ non-negative such that $\pk{ \epsilon_\delta(Y)>x}>0$  we have that
	$\mathcal{ B}_{ Z }^\delta (x) \in (0,\IF)$, hence the proof follows. \\
	{Assume now that  $\pk{S_0< \IF}\in (0,1)$.} In view of \nekorr{xhxh} we have
 \bqny{
 	\mathcal{ B}_{ Z }^ \delta(x)	&=& 	b\mathcal{ B}_{ V}^ \delta(x),
 }
 with $V(t)= Z(t)\lvert S_0 < \IF $, which is well-defined since $\pk{S_0<\IF}>0$ by the assumption.
 Since $S_0(V)< \IF$ almost surely and
 $ Y_*(t) =Y(t) \lvert S_0(\Theta )< \IF,t\inr$ by the proof above
 \bqny{
	\mathcal{ B}_{ Z }^ \delta(x)	&=& \pk{S_0< \IF}\E*{\frac{ \ind{ \epsilon_\delta (Y_*) > x}}{\epsilon_\delta (Y_* )}} \\
		 	&=& \E*{\frac{ \ind{ \epsilon_\delta (R \Theta) > x}}{\epsilon_\delta (R \Theta )} \ind{
		 			S_0(\Theta)< \IF}}.
}
In view of \cite[Lem 2.5, Cor 2.9]{PH2020} and \cite[Thm 3.8]{HBernulli}  and the above
\bqn{ \label{nukd}
	H_{ Z }^\delta =\mathcal{ B}_{ Z }^ \delta(0)=\E*{\frac{ 1}{\epsilon_\delta (R \Theta )}}=
	\E*{\frac{ 1}{\epsilon_\delta (R \Theta )} \ind{ S_0(\Theta)< \IF}}
}	
and thus  $\epsilon_\delta (R \Theta )< \IF$ implies $S_0< \IF$ almost surely. Hence the proof is complete.
\QED 	

\proofkorr{korrIm}
In view of \eqref{K1},
the representation \eqref{repBE} and the finiteness of $\mathcal{B}_Z^0(x)$ for all $x\ge 0$,
the monotone convergence theorem yields for all $x_0\ge 0$
$$ \lim_{x \downarrow x_0} \E*{\frac{\ind{ x_0 \le \epsilon_0 (Y )<x}}{\epsilon_0 (Y )} } =
\E*{\frac{\ind{\epsilon_0 (Y )=x_0}}{\epsilon_0 (Y )} }=0$$
consequently, since by  our assumption \Cref{L1} implies  $\pk{\epsilon_0 (Y ) \in (0,\IF)}=1$, then
$$\pk{\epsilon_0 (Y )=x_0 } = \E{ \ind{ \epsilon_0 (Y )=x_0}}=0$$
follows establishing the claim.
\QED

	\proofprop{ineq1}
	%Using \eqref{repBE} for all $x\ge 0$ we have
	%\bqny{
		%	\mathcal{B}_{ Z }^\delta (x) &=& \int_{x}^\infty  \frac{1}{y}dF(y)=\sum_{k=0}%^\infty\int_{x+k}^{x+k+1}  \frac{1}{y}dF(y),
		%}
	%$where $F$ is the distribution of $\epsilon_\delta(Y)$. Thus by $\frac{1}{x+k+1}%$<\frac{1}{y}\leq\frac{1}{x+k}$
	%for $y\in[x+k, x+k+1)$ we obtain the first claim. The second claim follows by %Markov's inequality.\\
	In order to prove \eqref{eqXhK} note first that
	%Applying Jensen inequality the following lower bound for Pickands constant follows
	%\begin{eqnarray}
	%	\mathcal{ B}_{ Z }^\delta (0)&\geq& \frac{1}{\E{\epsilon_\delta(Y)}}\nonumber\\
	%	&=&\frac{1}{\int_{\delta \Z^d}\int_1^\infty\Psi\left(\frac{\sigma(t)}{2}-\frac{\ln  r}{\sigma(t)}\right)r^{-2}\lambda_\delta(dt)dr}\,.\label{lowerbound}
	%\end{eqnarray}
	%%For $x>0$ that is Berman constant we can also use Jensen inequality as well.
	%Note that
	for any non-negative rv $U$ with df $G$ and $x\ge 0 $ such that $\pk{U>x}>0$
	\begin{eqnarray*}
		\frac{1}{\pk{U>x}}\int_x^\infty\frac{1}{y}dG(y)&\geq&\frac{\pk{U>x}}{\int_x^\infty ydG(y)}\ge \frac{\pk{U>x}}{\E{U}}.
	\end{eqnarray*}
	Consequently, we obtain for all $x>0$
	\begin{eqnarray*}
		\mathcal{ B}_{ Z }^\delta (x)&\geq& \frac{\pk{\epsilon_\delta(Y)>x}^2}{\E{\epsilon_\delta(Y)\mathbb{I}
				\{\epsilon_\delta(Y)>x\}}} \geq \frac{\pk{\epsilon_\delta(Y)>x}^2}{\E{\epsilon_\delta(Y)}}
		%		\\
		%	&=&\frac{\pk{\epsilon_\delta(Y)>x}^2}{\int_{\delta \Z^d}\int_1^\infty\Psi\left(\frac{\sigma(t)}{2}-\frac{\ln  r}{\sigma(t)}\right)r^{-2}\lambda_\delta(dt)dr}\,.
	\end{eqnarray*}
	establishing the proof
	of the lower bound \eqref{eqXhK}. The proof of the upper bound  follows from the fact that
	\[
	\mathcal{ B}_{ Z }^\delta (x)
	= \int_{x}^\infty  \frac{1}{y}dF(y)
	\leq
	\frac{1}{x}\int_{x}^\infty  dF(y)=x^{-}{\pk{\epsilon_\delta(Y)\geq x}},
	\]
	where $F$ is the distribution of $\epsilon_\delta(Y)$. This completes the proof.
	\QED
	
\proofprop{korrDe} Since $\mathcal{ B}_{ Z }^\delta(0)$ is the generalised Pickands constant $\mathcal{H}_{ Z }^\delta$, then the claim follows for $x=0$ from \cite{ZKE}.  In view of \eqref{paP}  we can assume  without loss of generality that $\pk{S_0 < \IF}=1$. Under this assumption, from the proof of \Cref{L1} we have that
$Y(t)\to 0$ almost surely as $\norm{t}\to \IF$. Hence for some $M$ sufficiently large $Y(t) < 1$ almost surely for all $t$ such that $\norm{t}> M$. Consequently, for all $\delta\ge 0$
$$\epsilon_\delta(Y) =  \int_{\delta Z ^d \cap [-M, M]^d} \ind{ Y(t)> 1} \lambda_\delta(dt).
$$
Moreover, $\epsilon_\delta(Y)< \IF$ almost surely for all $\delta \ge 0$ implying   $\epsilon_\delta(Y) \to \epsilon_0(Y)$ almost surely as $\delta \downarrow 0$. In view of \cite[Lem 2.5, Cor 2.9]{PH2020} and \cite[Thm 3.8]{HBernulli} for all  $\delta \ge 0$
$$ \mathcal{H}_{ Z }^\delta =\E{ 1/\epsilon_\delta(Y)  } .$$
Applying \cite[Thm 2]{ZKE}  and \eqref{nukd} yields
$$ \E{ 1/\epsilon_\delta(Y)  } = \mathcal{H}_{ Z }^{\delta  } \to  \mathcal{H}_{ Z }^0= \E{1/\epsilon_0(Y) }, \quad \delta  \downarrow 0.
$$
Hence $1/\epsilon_{\delta }(Y), \delta >0 $ is uniformly integrable and hence
$$\mathcal{B}_{ Z }^\delta(x)= \E*{ \frac{ \ind{\epsilon_\delta(Y) > x   }}{\epsilon_\delta(Y) }} \to
   \mathcal{B}_{ Z }^0(x), \quad \delta \downarrow 0  $$
   establishing  the proof.
\QED

\subsection{Proof of Theorem \ref{Th.rate}}\label{s.proof.rate}
Suppose that $V(t), t\in\R^d$ is a centered Gaussian field with stationary increments and variance function
$\sigma^2_V(\cdot)$ that satisfies {\bf A1-A2}.
Then, by stationarity of increments $\sigma^2_V(\cdot)$ is negative definite, which
combined with Schoenberg's theorem, implies that for each $u>0$
\[
R_u(s,t)\coloneqq\exp\left(-\frac{1}{2u^2}\sigma_V^2(s-t)\right), \ s,t\in\R^d
\]
is positive definite, and thus a valid covariance function of some centered \cEH{stationary}
Gaussian rf $X_u(t), t\in\R^d$, where $s-t$ is meant component-wise.
The proof of \Cref{Th.rate} is based on the analysis of the asymptotics of sojourn time of
$X_u(t)$. Since the idea of the proof is the same for continuous and discrete scenario,
in order to simplify notation, we consider next only the case $\delta=0$.

Before we proceed to the proof of \Cref{Th.rate}, we need the following lemmas, where
$Z(t)=\exp\left(V(t)-\frac{\sigma^2_V(t)}{2}\right)$ is as in \Cref{cern}, \Cref{MMS}.
\begin{lem}\label{l.sojourn}
For all $T>0$ and $x\ge0$
\begin{enumerate}[(i)]
	\item \label{thmK1}
\[
\lim_{u\to\infty}\frac{\pk*{  \int_{[0,T]^d  } \mathbb{I}(X_u(t)>u) dt>x   }}{\Psi(u)}=\mathcal{B}_Z([0,T]^d,x).
\]
\item \label{thmK2} For all $x\ge 0$
\[
\lim_{u\to\infty}\frac{\pk*{  \int_{[0,\ln (u)]^d  } \mathbb{I}(X_u(t)>u) dt>x   }}{(\ln (u))^d\Psi(u)}=\mathcal{B}_Z(x), \quad
\lim_{T\to\infty}\frac{\mathcal{B}_Z([0,T]^d,x)}{T^d}=\mathcal{B}_Z(x)\in(0,\infty).
\]
\end{enumerate}
\end{lem}
\prooflem{l.sojourn}
\Cref{thmK1} follows straightforwardly from \cite[Lem 4.1]{DHLM22}.
The proof of \Cref{thmK1} follows by the application of the double sum technique applied to the sojourn functional,
as demonstrated e.g., in \cite[Prop 3.1]{DHLM22}.
The claim in  \Cref{thmK2} follows by the same argument as its counterpart in \cite[Lem 4.2]{DHLM22}.
\QED

The following lemma is a slight modification of \cite[Lem 6.3]{Pit96} to the family $X_u,\ u>0$.
Let $\mathbf{i}=(i_1,...,i_d)$, with $i_1,...,i_d\in \{0,1,2,...\}$,
$\mathcal{R}_{\mathbf{i}}:=\prod_{k=1}^d [i_k T,(i_k+1)T]$
and
\begin{eqnarray*}
\widehat{\mathcal{K}}&:=&\{\mathbf{i}=(i_1,...,i_d): 0\le i_k, (i_k-1) T\le \ln (u), k=1,...,d    \},\\
\widecheck{\mathcal{K}}&:=&\{\mathbf{i}=(i_1,...,i_d): 0\le i_k T\le \ln (u), k=1,...,d    \}.
\end{eqnarray*}

\begin{lem}\label{l.double}
There exists a constant $C\in(0,\infty)$ such that for sufficiently large $u$,
for all $\mathbf{i},\mathbf{j}\in\widehat{\mathcal{K}}, \mathbf{i}\neq\mathbf{j}$ we have
\[
\pk*{\max_{t\in \mathcal{R}_{\mathbf{i}}} X_u(t)>u,\max_{t\in \mathcal{R}_{\mathbf{j}}} X_u(t)>u}
\le
C T^{2d}\exp\left(  -\frac{1}{8} \inf_{t\in \mathcal{R}_{\mathbf{i}}, s\in \mathcal{R}_{\mathbf{j}}}
\sigma_V^2(t-s) \right) \Psi(u).
\]
\end{lem}
{\it Proof of Theorem \ref{Th.rate}.}
The proof consists of two steps, where we find an asymptotic upper and lower bound for the ratio
\[
\frac{\pk*{  \int_{[0,\ln (u)]^d  } \mathbb{I}(X_u(t)>u) dt>x   }}{(\ln (u))^d\Psi(u)},
\]
as $u\to\infty$. We note that by \Cref{l.sojourn} the limit, as $u\to\infty$, of the above fraction
is positive and finite.  \\
{\it \underline{Asymptotic upper bound.}}
If  $T>0$, then for sufficiently large $u$
\begin{eqnarray}
\lefteqn{\pk*{  \int_{[0,\ln (u)]^d  } \mathbb{I}(X_u(t)>u) dt>x   }}\nonumber\\
&\le&
\pk*{ \sum_{\mathbf{i}\in \widehat{\mathcal{K}}} \int_{\mathcal{R}_{\mathbf{i}}  } \mathbb{I}(X_u(t)>u) dt>x   }\nonumber\\
&\le&
\pk*{ \exists_{\mathbf{i}\in \widehat{\mathcal{K}}} \int_{\mathcal{R}_{\mathbf{i}}  } \mathbb{I}(X_u(t)>u) dt>x   }
+
\pk*{ \exists_{\mathbf{i},\mathbf{j}\in \widehat{\mathcal{K}},\mathbf{i}\neq\mathbf{j} }
\max_{t\in\mathcal{R}_{\mathbf{i}}  } X_u(t)>u, \max_{t\in\mathcal{R}_{\mathbf{j}}  } X_u(t)>u }\nonumber\\
&\le&
\sum_{\mathbf{i}\in \widehat{\mathcal{K}}}
\pk*{  \int_{\mathcal{R}_{\mathbf{i}}  } \mathbb{I}(X_u(t)>u) dt>x   }
+
\pk*{ \exists_{\mathbf{i},\mathbf{j}\in \widehat{\mathcal{K}},\mathbf{i}\neq\mathbf{j} }
\max_{t\in\mathcal{R}_{\mathbf{i}}  } X_u(t)>u, \max_{t\in\mathcal{R}_{\mathbf{j}}  } X_u(t)>u }\nonumber\\
&\le&
\left\lceil \frac{(\ln (u))^d}{T^d}\right\rceil
\pk*{  \int_{[0,T]^d  } \mathbb{I}(X_u(t)>u) dt>x   }
+
\pk*{ \exists_{\mathbf{i},\mathbf{j}\in \widehat{\mathcal{K}},\mathbf{i}\neq\mathbf{j} }
\max_{t\in\mathcal{R}_{\mathbf{i}}  } X_u(t)>u, \max_{t\in\mathcal{R}_{\mathbf{j}}  } X_u(t)>u },\nonumber\\
\label{up.1}
\end{eqnarray}
where $\lceil\cdot\rceil$ is the ceiling function and the last inequality above follows from the \cEH{stationarity} of $X_u$. Using again the \cEH{stationary} of $X_u$, we obtain
\begin{eqnarray}
\lefteqn{
\pk*{ \exists_{\mathbf{i},\mathbf{j}\in \widehat{\mathcal{K}},\mathbf{i}\neq\mathbf{j} }
\max_{t\in\mathcal{R}_{\mathbf{i}}  } X_u(t)>u, \max_{t\in\mathcal{R}_{\mathbf{j}}  } X_u(t)>u }}\nonumber\\
&\le&
\sum_{\mathbf{i},\mathbf{j}\in \widehat{\mathcal{K}},\mathbf{i}\neq\mathbf{j} }\pk*{
\max_{t\in\mathcal{R}_{\mathbf{i}}  } X_u(t)>u, \max_{t\in\mathcal{R}_{\mathbf{j}}  } X_u(t)>u }\label{d.1}\\
&\le&
\left\lceil \frac{(\ln (u))^d}{T^d}\right\rceil
\sum_{\mathbf{k}\in \widehat{\mathcal{K}},\mathbf{k}\neq\mathbf{0}  }\pk*{
\max_{t\in\mathcal{R}_{\mathbf{0}}  } X_u(t)>u, \max_{t\in\mathcal{R}_{\mathbf{k}}  } X_u(t)>u }\nonumber\\
&=&
\left\lceil \frac{(\ln (u))^d}{T^d}\right\rceil
\left(
\sum_{\mathbf{k}\in \widehat{\mathcal{K}},\mathbf{k}\neq\mathbf{0}, \mathcal{R}_{\mathbf{0}}\cap \mathcal{R}_{\mathbf{k}}\neq\emptyset  }\pk*{
\max_{t\in\mathcal{R}_{\mathbf{0}}  } X_u(t)>u, \max_{t\in\mathcal{R}_{\mathbf{k}}  } X_u(t)>u }\right. \nonumber\\
&&+
\left.
\sum_{\mathbf{k}\in \widehat{\mathcal{K}},\mathbf{k}\neq\mathbf{0}, \mathcal{R}_{\mathbf{0}}\cap \mathcal{R}_{\mathbf{k}}=\emptyset  }
\pk*{ \max_{t\in\mathcal{R}_{\mathbf{0}}  } X_u(t)>u, \max_{t\in\mathcal{R}_{\mathbf{k}}  } X_u(t)>u }
\right)\nonumber\\
&=:&\left\lceil \frac{\ln ^d(u)}{T^d}\right\rceil
\left(\Sigma_1+\Sigma_2\right). \label{up.2}
\end{eqnarray}
Next, by \Cref{l.double}, for sufficiently large $T,u$ and some $\rm{Const_0}>0$
\begin{eqnarray}
\Sigma_2 &\le&
C T^{2d}
\sum_{\mathbf{k}\in \widehat{\mathcal{K}},\mathbf{k}\neq\mathbf{0}, \mathcal{R}_{\mathbf{0}}\cap \mathcal{R}_{\mathbf{k}}=\emptyset  }
\exp\left( -\frac{1}{8}\sigma_V^2(T\mathbf{k}) \right) \Psi(u)\nonumber\\
&\le&
C T^{2d}
\sum_{\mathbf{k}>\mathbf{0}  }
\exp\left( -\rm{Const_0} T^{\alpha_\infty}\lVert\mathbf{k}\rVert^{\alpha_\infty} \right)\Psi(u)\nonumber\\
&\le&
{\rm{Const_1}} T^{2d} \exp\left( - T^{\alpha_\infty/2} \right)\Psi(u).\label{up.3}
\end{eqnarray}

The upper bound for $\Sigma_1$ follows by a similar argument as used in the proof of
\cite[Lem 6.3]{Pit96}, thus we explain only main steps of the argument.
For a while, consider the following probability
\[\pk*{ \max_{t\in\mathcal{R}_{\mathbf{0}}  } X_u(t)>u, \max_{t\in\mathcal{R}_{{(1,0,...,0)}}  } X_u(t)>u }.\]
Then, for each $\varepsilon>0$ and sufficiently large $T$, $u$,
\begin{eqnarray}
\lefteqn{\pk*{ \max_{t\in\mathcal{R}_{\mathbf{0}}  } X_u(t)>u, \max_{t\in\mathcal{R}_{(1,0,...,0)}  } X_u(t)>u }}\nonumber\\
&\le&
\pk*{ \max_{t\in [0,T^\varepsilon]\times [0,T]^{d-1}} X_u(t)>u }
+
\pk*{ \max_{t\in [0,T]^{d}  } X_u(t)>u, \max_{t\in [T^\varepsilon,T^\varepsilon+T] \times [0,T]^{d-1} } X_u(t)>u }\nonumber\\
&\le&
{\rm Const_2} T^{d-1+\varepsilon}\Psi(u)
+
{\rm Const_3} T^{2d}\exp\left( -T^{\varepsilon/2}\right),\label{up.4}
\end{eqnarray}
where the above inequality follows by \Cref{l.double} and
\begin{eqnarray*}
%\lefteqn{
\lim_{u\to\infty}
\frac{\pk*{ \max_{t\in [0,T^\varepsilon]\times [0,T]^{d-1}} X_u(t)>u }}
     {\Psi(u)}
&\le&
\lceil T \rceil^{(d-1)(1-\epsilon)}
\lim_{u\to\infty}
 \frac{\pk*{ \max_{t\in [0,T^\varepsilon]^d} X_u(t)>u }}
     {\Psi(u)}\\
     &=&\lceil T \rceil^{(d-1)(1-\epsilon)} \mathcal{B}_Z([0,T^\varepsilon]^d,0)\\
     &\le& {\rm Const_4} T^{d-1+\varepsilon},
\end{eqnarray*}
which  is a consequence of the \cEH{stationarity} of $X_u$ and  statement (i) of \Cref{l.sojourn} applied to $x=0$.
Again, by the \cEH{stationarity} of $X_u$ we can obtain the  bound as  in (\ref{up.4})
uniformly for all the summands in $\Sigma_1$.

Application of the bounds (\ref{up.2}), (\ref{up.3}), (\ref{up.4}) to (\ref{up.1}) leads to the following upper estimate
\begin{eqnarray}
\lefteqn{
\limsup_{u\to\infty}
\frac{\pk*{  \int_{[0,\ln (u)]^d  } \mathbb{I}(X_u(t)>u) dt>x   }}{\ln ^d(u)\Psi(u)}
}\nonumber\\
&\le&
\frac{\mathcal{B}_Z([0,T]^{d},x)}{T^d}+
{\rm Const_4}
\frac{1}{T^d}\left(
T^{d-1+\varepsilon}
+
T^{2d} \exp\left( - T^{\alpha_\infty/2} \right)
+
T^{2d}\exp\left( -T^{\varepsilon/2}\right)
\right),\label{sup.1}
\end{eqnarray}
which is valid for all  $\varepsilon>0$ and $T$ sufficiently large.\\
{\it \underline{Asymptotic lower bound.}} Taking  $T>0$, for sufficiently large $u$
\begin{eqnarray}
\lefteqn{\pk*{  \int_{[0,\ln (u)]^d  } \mathbb{I}(X_u(t)>u) dt>x   }}\nonumber\\
&&\ge
\pk*{ \sum_{\mathbf{i}\in \widecheck{\mathcal{K}}} \int_{\mathcal{R}_{\mathbf{i}}  } \mathbb{I}(X_u(t)>u) dt>x   }\nonumber\\
&&\ge
\pk*{ \exists_{\mathbf{i}\in \widecheck{\mathcal{K}}} \int_{\mathcal{R}_{\mathbf{i}}  } \mathbb{I}(X_u(t)>u) dt>x   }\nonumber\\
%-
%\pk*{ \exists_{\mathbf{i},\mathbf{j}\in \widecheck{\mathcal{K}},\mathbf{i}\neq\mathbf{j} }
%\max_{t\in\mathcal{R}_{\mathbf{i}}  } X_u(t)>u, \max_{t\in\mathcal{R}_{\mathbf{j}}  } X_u(t)>u }\\
&&\ge
\sum_{\mathbf{i}\in \widecheck{\mathcal{K}}}
\pk*{  \int_{\mathcal{R}_{\mathbf{i}}  } \mathbb{I}(X_u(t)>u) dt>x   }
-
\sum_{\mathbf{i},\mathbf{j}\in \widecheck{\mathcal{K}},\mathbf{i}\neq\mathbf{j} }
\pk*{
\max_{t\in\mathcal{R}_{\mathbf{i}}  } X_u(t)>u, \max_{t\in\mathcal{R}_{\mathbf{j}}  } X_u(t)>u }\label{bon1}\\
&&\ge
\left\lfloor \frac{\ln ^d(u)}{T^d}\right\rfloor
\pk*{  \int_{[0,T]^d  } \mathbb{I}(X_u(t)>u) dt>x   }
-
\sum_{\mathbf{i},\mathbf{j}\in \widecheck{\mathcal{K}},\mathbf{i}\neq\mathbf{j} }
\pk*{
\max_{t\in\mathcal{R}_{\mathbf{i}}  } X_u(t)>u, \max_{t\in\mathcal{R}_{\mathbf{j}}  } X_u(t)>u }\,,\nonumber\\
\label{down.1}
\end{eqnarray}
where in (\ref{bon1}) we used Bonferroni inequality.

Using that $\widecheck{\mathcal{K}}\subset \widehat{\mathcal{K}}$ with  the upper bound for
\[
\sum_{\mathbf{i},\mathbf{j}\in \widehat{\mathcal{K}},\mathbf{i}\neq\mathbf{j} }
\pk*{
\max_{t\in\mathcal{R}_{\mathbf{i}}  } X_u(t)>u, \max_{t\in\mathcal{R}_{\mathbf{j}}  } X_u(t)>u }
\]
derived in (\ref{d.1}), we conclude that
for each $T$
sufficiently large and $\varepsilon>0$,
\begin{eqnarray}
\lefteqn{
\liminf_{u\to\infty}
\frac{\pk*{  \int_{[0,\ln (u)]^d  } \mathbb{I}(X_u(t)>u) dt>x   }}{(\ln (u))^d\Psi(u)}}\nonumber\\
&\ge&
\frac{\mathcal{B}_Z([0,T]^{d},x)}{T^d}-
{\rm Const_4}
\frac{1}{T^d}\left(
T^{d-1+\varepsilon}
+
T^{2d} \exp\left( - T^{\alpha_\infty/2} \right)
+
T^{2d}\exp\left( -T^{\varepsilon/2}\right)
\right). \label{inf.1}
\end{eqnarray}

Thus, by statement (ii) of \Cref{l.sojourn} combined with (\ref{sup.1}) and (\ref{inf.1}),
in view of the fact that $\varepsilon$ can take any value in $(0,1)$, we arrive at
\begin{eqnarray*}
\lim_{T\to\infty}\left|\mathcal{ B}_{ Z }(x)-\frac{\mathcal{ B}_{ Z }([0,T]^d,x)}{T^d}\right|
T^\lambda=0
\end{eqnarray*}
for all $\lambda\in (0,1)$ establishing the proof.
\QED

\proofprop{prop.bound}
The idea of the proof is to analyze the asymptotic upper and lower bound of
\[
\pk{\epsilon_\delta(Y)>x}
\]
as $x\to\infty$ and then to apply \Cref{ineq1}.
In order to simplify the notation, we consider only the case $\delta=0$.
Let $Z(t)=V(t)-\sigma^2_V(t)/2, t\in \R$ with $V$ a centered Gaussian process
with stationary increments that satisfies {\bf A1-A2} and $\mathcal{W}$ an independent of $V$ exponentially
distributed rv with parameter $1$.
\\
{\it \underline{Logarithmic upper bound.}}
Let $A\in (0,1/2)$.
We begin with an observation that
\BQN \lefteqn{ \pk*{\epsilon_\delta(Y)>x}}\nonumber\\
&=&
\pk*{ \int_{\R}\mathbb{I}\{ \mathcal{W}+ V(t)-\sigma^2_V(t)/2>0 \} dt >x}\nonumber\\
&\le&
\pk*{\mathcal{W}\le A\sigma^2_V(x/2), \int_{\R}\mathbb{I}\{ A\sigma^2_V(x/2)+ V(t)-\sigma^2_V(t)/2>0 \} dt >x}\nonumber\\
&& +\pk*{\mathcal{W}>A\sigma^2_V(x/2)}\nonumber\\
&\le&
e^{-A\sigma^2_V(x/2)}
+
\pk*{  \int_{\R}\mathbb{I}\{ A\sigma^2_V(x/2)+ V(t)-\sigma^2_V(t)/2>0 \} dt >x}\nonumber\\
&\le&
e^{-A\sigma^2_V(x/2)}
+
\pk*{  \sup_{t\in\left(-\infty,-x/2]\cup [x/2,\infty\right)} V(t)-\sigma^2_V(t)/2>  -A\sigma^2_V(x/2)}\nonumber\\
&\le&
e^{-A\sigma^2_V(x/2)}
+
2\pk*{  \exists_{t\in \left[x/2,\infty\right)} V(t)>\left(\frac{1}{2}  -A\right)\sigma^2_V(t)}\label{line}\\
&=&
e^{-A\sigma^2_V(x/2)}
+
2\pk*{  \exists_{t\in \left[x/2,\infty\right)} \frac{V(t)}{\sigma_V^2(t)}>
\left(\frac{1}{2}  -A\right)}\, ,\label{upp.2}
\EQN
where in (\ref{line}) we used that $\{V(-t),t\ge 0\}\stackrel{d}{=} \{V(t),t\ge0\} $
and the assumption that $\sigma^2_V$ is increasing.
%combined with $A<1/2$ ensures that
%for sufficiently large $x$ and each $t\ge x/2$, we have
%$\sigma^2_V(t)/2-A\sigma^2_V(x/2)>0$.
%Moreover, one can check that
%\[
%\max_{t\in \left[x/2,\infty\right)} Var \left( \frac{V(t)}{\sigma^2_V(t)/2  -A\sigma^2_V(x/2)}\right)  =
%\frac{1}{\left(\frac{1}{2}-A\right)^2}\frac{1}{\sigma_V^2(x/2)}
%\]
%as $x\to\infty$.
Next, by {\bf A1}, for sufficiently large $x$ and $s,t\ge  x/2$ such that $|t-s|\le 1$
$$Cov\left(\frac{V(t)}{\sigma_V(t)},\frac{V(s)}{\sigma_V(s)}\right) \ge
\exp\left(-|t-s|^{\alpha_0/2}\right)=:Cov\left(Z(t),Z(s)\right),$$
where $Z$ is some centered stationary Gaussian process.
Hence, by Slepian inequality (see, e.g., Corollary 2.4 in \cite{Adl90})
\begin{eqnarray*}
\pk*{  \exists_{t\in \left[x/2,\infty\right)} \frac{V(t)}{\sigma^2_V(t)}>\left(\frac{1}{2}  -A\right)}
&\le&
\sum_{k=0}^\infty
\pk*{  \exists_{t\in [x/2+k,x/2+k+1]} \frac{V(t)}{\sigma^2_V(t)}>
\left(\frac{1}{2}  -A\right)}\\
&\le&
\sum_{k=0}^\infty
\pk*{  \exists_{t\in [0,1]} Z(t)>
\left(\frac{1}{2}  -A\right){\sigma_V(x/2+k)}}
\end{eqnarray*}
and
by Landau-Shepp (see, e.g.,  \cite[Eq.\ (2.3)]{Adl90}), uniformly with respect to $k$
\[
\lim_{x\to\infty}\frac{\ln \left(\pk*{  \exists_{t\in [0,1]} Z(t)>
\left(\frac{1}{2}  -A\right){\sigma_V(x/2+k)}}
\right)}
{\sigma_V^2(x/2+k)}=-\frac{1}{2}\left(\frac{1}{2}-A\right)^2\,.
\]
The above implies that
\[
\lim_{x\to\infty}\frac{\ln
\left(\pk*{  \exists_{t\in \left[x/2,\infty\right)} \frac{V(t)}{\sigma^2_V(t)}>\left(\frac{1}{2}  -A\right)}
\right)}
{\sigma_V^2(x/2)}\le
-\frac{1}{2}\left(\frac{1}{2}-A\right)^2\ .
\]
Thus, in order to optimize the value of $A$ in (\ref{upp.2})  it suffices now to solve
\[
\left(\frac{1}{2}-A\right)^2=2A
\]
that leads to (recall that $A<1/2$) 
\[
A=\frac{3-2\sqrt{2}}{2}.
\]
Hence
\[
\lim_{x\to\infty}\frac{\ln \left(\pk{\epsilon_\delta(Y)>x}\right)}{\sigma_V^2(x/2)}\le - \frac{3-2\sqrt{2}}{2},
\]
which combined with (\ref{eqXhK}) in \Cref{ineq1} completes the proof of the logarithmic upper bound.\\
{\it \underline{Logarithmic lower bound.}}
Taking  $A> 1/2$ we have 
\begin{eqnarray*}
\pk{\epsilon_\delta(Y)>x}
&=&
\pk*{ \int_{\R}\mathbb{I}\{ \mathcal{W}+ V(t)-\sigma^2_V(t)/2>0 \} dt >x}\\
&\ge&
\pk*{\mathcal{W}>A\sigma^2_V(x/2), \int_{\R}\mathbb{I}\{ A\sigma^2_V(x/2)+ V(t)-\sigma^2_V(t)/2>0 \} dt >x}\\
&\ge&
\pk*{\mathcal{W}>A\sigma^2_V(x/2)}
\pk*{ \inf_{t\in [-x/2, x/2]} V(t)> -(A-1/2)\sigma^2_V(x/2) }\\
&=&
e^{-A\sigma^2_V(x/2)}
\left(1-\pk{ \sup_{t\in [-x/2, x/2]} V(t)> (A-1/2)\sigma^2_V(x/2) }\right).
\end{eqnarray*}
Using that
\begin{eqnarray}
\lefteqn{\pk*{ \sup_{t\in [-x/2, x/2]} V(t)> (A-1/2)\sigma^2_V(x/2) }}\nonumber\\
&\le&
2\sum_{i \in\{0,...,\lfloor x/2 \rfloor-1\}}
\pk*{ \sup_{t\in [i,i+1]} V(t)> (A-1/2)\sigma^2_V(x) }\\
&& +
2\pk*{ \sup_{t\in [\lfloor \zm{x/2} \rfloor,\zm{x/2}]} V(t)> (A-1/2)\sigma^2_V(x) }\label{sum}
\end{eqnarray}
and the fact that by the stationarity of increments of $V$
\[
\E*{\sup_{t\in [i,i+1]} V(t)}=
\E*{\sup_{t\in [i,i+1]} (V(t)-V(i))+V(i)}
=
\E*{\sup_{t\in  [0,1]} V(t)}=:\mu< \IF
\]
we can apply Borell inequality (e.g.,  \cite[Thm 2.1]{Adl90}) uniformly for all the summands in (\ref{sum})
to get
that for sufficiently large $x$ (recall that $\sigma^{2}_V$ is supposed to be increasing)
\begin{eqnarray*}
\pk*{ \sup_{t\in [-x/2, x/2]} V(t)> (A-1/2)\sigma^2_V(x/2) }
&\le&
4 (x+1)\exp\left(-\frac{((A-1/2)\sigma^{2}_V(x/2)\zm{-}\mu)^2}{2 \sigma^2_V(x/2)}     \right)\\
&\le&\exp\left(- \frac{(A-1/2)^2\sigma^{2}_V(x/2)}{4}\right)\to 0
\end{eqnarray*}
as $x\to\infty$.
%\footnote{\zm{$\sigma^2(t)$ is less than 2 on every subsquare? What happens with $\E{\sup V(t)}$? Ok we can use
%Borel-Piterbarg inequality but with $a>0$ and $a$ should depend on $x$}\\
%K: I have added more explanation, should be ok.}
Hence we
arrive at
\[
\liminf_{x\to\infty}\frac{\ln (\pk{\epsilon_\delta(Y)>x})}{\sigma^2_V(x/2)}\ge -A,
\]
which combined with (\ref{eqXhK}) in \Cref{ineq1} and the fact that, by the proof of the logarithmic upper bound
$\E{\epsilon_\delta(Y)}<\infty$  implies
%\footnote{\zm{Do we not need to assume that $\E{\epsilon_\delta(Y)}<\infty$? Ok it is finite by A2 and the formula.}\\
%K: Yes, finiteness follows from the upper bound. I will change the order of the proofs. I hope, now it should be clearer.}
\[
\liminf_{x\to\infty}\frac{\ln (\mathcal{ B}_{ Z } (x))}{\sigma^2_V(x)}\ge -2A
\]
for all $A>1/2$. This completes the proof.
\QED

\COM{
\prooflem{lemesx}
Using inequality (\ref{eqXhK}) we have $\pk{\epsilon_\delta(Y)>x}\leq (\mathcal{ B}_{ Z }^\delta (x)\E{\epsilon_\delta(Y)})^{1/2}$. Thus for $s>0$
\begin{eqnarray*}
\E{e^{s\epsilon_\delta(Y)}}&=&s\int_{-\infty}^\infty e^{sx}\pk{\epsilon_\delta(Y)>x}dx\\
&=&1+s\int_{0}^\infty e^{sx}\pk{\epsilon_\delta(Y)>x}dx\\
&\leq&1+ s(\E{\epsilon_\delta(Y)})^{1/2}\int_0^\infty e^{sx}(\mathcal{ B}_{ Z }^\delta (x))^{1/2}dx\,.
\end{eqnarray*}
The proof is complete.
\QED
}

\bibliographystyle{ieeetr}
\bibliography{EB}
\end{document}